\documentclass[12pt, twoside]{article}
\usepackage{amsmath,amsthm,amssymb}
\usepackage{times}
\usepackage{enumerate}
\usepackage{extarrows}
\usepackage{color}
\usepackage{tabularx}
\usepackage{multirow}
\usepackage{graphicx}

\def\titlerunning#1{\gdef\titrun{#1}}
\makeatletter
\def\author#1{\gdef\autrun{\def\and{\unskip, }#1}\gdef\@author{#1}}
\def\address#1{{\def\and{\\\hspace*{18pt}}\renewcommand{\thefootnote}{}%
\footnote {#1}}%
\markboth{\autrun}{\titrun}}
\makeatother
\def\email#1{e-mail: #1}

\def\keywords#1{\par\medskip
\noindent\textbf{Keywords.} #1}

\newtheorem{theorem}{Theorem}[section]
\newtheorem{lemma}{Lemma}[section]

\newtheorem{remark}{Remark}[section]
\newtheorem{corollary}{Corollary}[section]

\newcommand{\Proof}{\begin{proof}}
\newcommand{\End}{\end{proof}}

\numberwithin{equation}{section}

\newcommand{\PreserveBackslash}[1]{\let\temp=\\#1\let\\=\temp}
\newcolumntype{C}[1]{>{\PreserveBackslash\centering}p{#1}}
\newcolumntype{R}[1]{>{\PreserveBackslash\raggedleft}p{#1}}
\newcolumntype{L}[1]{>{\PreserveBackslash\raggedright}p{#1}}

\newcolumntype{I}{!{\vrule width 1pt}}
\newlength\savedwidth

\begin{document}

%%%%% To ease editing, add:

\baselineskip=15pt

%%%%%%%%%%%%%%%%

%% In the running head, give an abbreviation of the title.
\titlerunning{ }

\title{The elliptical invariant tori of nearly integrable Hamiltonian system through symplectic algorithms}

\author{Zaijiu Shang ,  Yang Xu}

\date{\today}

\maketitle

\address{Zaijiu Shang:  Institute of Academy of Mathematics and Systems Science, Chinese Academy of Sciences, Beijing 100080, China; \email{zaijiu@amss.ac.cn}
		\and Yang Xu: School of Mathematical Sciences, Fudan University, Shanghai 200433, China; \email{xuyang$\_$@fudan.edu.cn}}

\begin{abstract}
In this paper we apply symplectic algorithms to nearly integrable Hamiltonian system, and prove it can maintain lots of elliptic lower dimensional invariant tori. 
 We are committed to consider the elliptic lower dimensional invariant tori for symplectic mapping with a small twist under the R\"{u}ssmann's non-degenerate condition, and focus on its measure estimation. And then apply it to the nearly integrable Hamiltonian system to obtain lots of elliptic lower dimensional invariant tori.
 
%    We are concerned with the stability of viscosity solutions to contact Hamilton-Jacobi equation 
%    \begin{equation*}
%	 H(x, \partial_x u ,u)=0,\quad x\in \mathbb{S} ,
%     \end{equation*}
%     where $\mathbb{S}$ is the unit circle and $H$ is of class $C^\infty$ and satisfies Tonelli conditions with respect to the argument $p$.
%     By analyzing the asymptotic behavior of the Lax-Oleinik dynamical system $(C(\ci ,\R),\{T^-_t\}_{t\geqslant 0})$, we find that the stability of solutions depends with the sign of constant $\mu$,  which is the weighted average of $ \partial H/\partial u$  along the circle $\mathbb{S}$. Moreover, we estimate the rate of convergence of asymptotically stable for $\mu >0$.
%
%At last, we also study the existence and multiplicity of nontrivial time periodic viscosity solutions  without monotonicity assumption.
\end{abstract}
	\keywords{nearly integrable, Symplectic algorithm, elliptic lower dimensional invariant tori, R\"{u}ssmann's non-degeneracy}

\tableofcontents

\section{Introduction}
\label{se1}

In section \ref{se1}, we briefly describe three  systems, they are the symplectic mapping system with a small twist (\ref{e01}), the nearly integrable Hamiltonian system (\ref{003}) and its corresponding symplectic difference scheme (\ref{004}). In section \ref{se2}, 
we give some notation. Section \ref{se3} contains the main conclusions. As for  the symplectic mapping system with a small twist (\ref{e01}), section \ref{se4} and section \ref{me1} are proof of the existence of elliptic lower dimensional invariant tori and measure estimation respectively. And applying this to the symplectic difference scheme is covered in section \ref{se6}.

Consider the symplectic mapping with a small twist
\begin{equation*}
    F :(x,u,y,v) \in \mathbb{T}^n \times W \times V \times W 
    \to (\hat{x},\hat{u},\hat{y},\hat{v}) \in
    \mathbb{T}^n \times \mathbb{R}^m \times  \mathbb{R}^n \times  \mathbb{R}^m
\end{equation*}
where $\mathbb{T}^n$ is the usual $n$-torus, $V$ is a   closed bounded and connected domain, $V \subseteq \mathbb{R}^n$ is a domain around the origin. The standard symplectic form is $\sum  \limits_{i=1}^n d x_i \land d y_i + \sum  \limits_{j=1}^m d u_j \land d v_j$. 
$F$ can be expressed implicitly as 
\begin{align} \label{e01}
      \hat{x} & =x+ \partial_{\hat{y}}t H(x,u,\hat{y},\hat{v})     & 
    y & = \hat{y} + \partial_{x}t H(x,u,\hat{y},\hat{v}) \notag \\
    \hat{u} & =u+ \partial_{\hat{v}}t H(x,u,\hat{y},\hat{v})      & 
    v & = \hat{v} + \partial_{u}t H(x,u,\hat{y},\hat{v})
\end{align}
where $t H = t N + t P$, $t$ is the time step, and $t P$ is viewed as a perturbation. The normal form $t N$ is expressed as
\begin{equation} \label{002}
    t N(x,u,\hat{y},\hat{v})= t h(\hat{y}) + \langle A_t u,\hat{v} \rangle +  \frac{1}{2} \langle B_t u,u \rangle + \frac{1}{2} \langle C_t \hat{v},\hat{v} \rangle
\end{equation}
For fixed $t$, $A_t = t(A-I), B_t = t B, C_t = t C $ are constant matrices, where $I$ is the identity matrix of order $m$. Suppose $B$ and $C$ are symmetric.

More importantly, we can apply it to nearly integrable Hamiltonian system:
\begin{equation}  \label{003}
    F^{'} : \dot{p} =  \frac{\partial H_\epsilon} {\partial q} (p,q)  , \   \dot{q} = - \frac{\partial H_\epsilon}{\partial p} (p,q),  \   (p,q) \in D
\end{equation}
where $H_\epsilon(p,q) = H_0 (p,q) + \epsilon H_1(p,q,\epsilon)$ ,
$H_0(x,u,y,v)=  h(y) + \langle A u,v \rangle +  \frac{1}{2} \langle B u,u \rangle + \frac{1}{2} \langle C v,v \rangle$, $p=(x,u)$, $q=(y,v)$.
$D \in {\mathbb{R}}^{2n}$ is a bounded  connected open domain, the dot is the derivative with respect to time $t$. 
Without loss of generality, we can assume $ \lvert H_1  \rvert _ D \leq M_1$, where $\lvert \cdot \lvert_ D$ denotes the supremum norm on $D$, $M_1$ is a constant, and we can ignore the parameter expression of $\epsilon$. 

Let us focus on the symplectic algorithm which is compatible with (\ref{003}). By Lemma 3.1 in \cite{r25}, the system $G_{H_ \epsilon}^t : (p,q) \to (\hat{p}, \hat{q}) $ can be expressed as follows:
\begin{equation} \label{004}
    G_{H_ \epsilon}^t :
    \begin{cases}
    \hat{p} = p + t  \frac{\partial   H_0} {\partial q} (p, \hat{q})
    + t \epsilon \frac{\partial   H_1} {\partial q} (p, \hat{q}) 
    + t^{s+1} \frac{ \partial{P_1}}{\partial q}(p, \hat{q}) \\
    \hat{q} = q - t  \frac{\partial   H_0} {\partial \hat{p}} (p, \hat{q}) 
    - t \epsilon \frac{\partial   H_1} {\partial \hat{p}} (p, \hat{q})
    - t^{s+1} \frac{ \partial{P_1}}{\partial \hat{p}}(p, \hat{q})
    \end{cases}
\end{equation}
for $(p, \hat{q}) \in D_{\frac{r}{4},\frac{s}{4}}$, where $s$ is a positive constant. 
Here $\delta$ is a sufficiently small positive number such that $P_1$ , which depends on time step $t$, is well defined and real analytic on $D_{\frac{r}{4},\frac{s}{4}}$ for $t \in [0,\delta]$, where $D_{r,s}$ is a complex neighborhood of $D$ with
$
    D_{r,s} = \{p : \lvert p \rvert < r \} \times \{q : \lvert \Im q \rvert < s \}
    \subset \mathbb{C}^n \times \mathbb{C}^n
$. 
By Lemma 3.1 and Lemma 3.3 in \cite{r25},  
$ \big \lvert \frac{\partial P_1}{\partial {\hat{p}}}  \big \rvert _ {\frac{r}{4},\frac{s}{4}} \leq M_2 $, $\big \lvert \frac{\partial P_1}{\partial q} \big \rvert _ {\frac{r}{4},\frac{s}{4}} \leq M_2 $, where $\lvert \cdot \lvert_ {\frac{r}{4},\frac{s}{4}}$ denotes the supremum norm on $D_{\frac{r}{4},\frac{s}{4}}$, $M_2$ is a constant. 
And similar to  \cite{r25}, fix $(p_0, \hat{q_0}) \in D$, let $P_1(p_0, \hat{q_0}) =0$, for 
$(p, \hat{q}) \in D_{\frac{r}{4},\frac{s}{4}}$,
we have $\lvert  P_1 \rvert _ {\frac{r}{4},\frac{s}{4}} \leq 2 n M_2 l_* $ for $t \in [0,\delta]$, where $l_*$ is an upper bound of the length of the shortest curves from $(p_0, \hat{q_0})$ to $(p, \hat{q})$ in $D_{\frac{r}{4},\frac{s}{4}}$.

By the way, for the nearly integrable Hamiltonian system \eqref{003}, there exists a generating function $\widetilde{S}$ such that  $H_1(p,q)  - \widetilde{S}(p,q) \sim o(t^s)$.

First we research the symplectic mapping system with a small twist (\ref{e01}), find its elliptic lower dimensional invariant tori, and give measure estimations. And then we applied it to the nearly integrable Hamiltonian system (\ref{003}). For the symplectic difference scheme  (\ref{004}) which is compatible with (\ref{003}), $S_0$ is regarded as $N$, $\epsilon H_1 + t^s P_1$ as $P$, and the conclusion of system (\ref{e01}) is applied. So we pay more attention to the symplectic mapping system with a small twist (\ref{e01}).

\section {Preprocessing}
\label{se2}
As for the symplectic mapping system with a small twist (\ref{e01}), if $P=0$, then
\begin{align*}
    \hat{x} & =x+ t\omega(\hat{y})    
   & y & = \hat{y}    \\
    \hat{u} & = \underline{A_t} u+ C_t \hat{v} 
   & v & = \underline{A_t}^T \hat{v} + B_t u
\end{align*}
where $\omega (\hat{y}) = h^{'} (\hat{y}) $, $\underline{A_t} = I+ A_t = I+t(A-I)$, $\underline{A_t}^T$ is the transpose of $\underline{A_t}$.

For $\underline{A_t}$ is non-singular, $F$ can be expressed as 
\begin{align*}
    \hat{x} & =x+ t\omega(y)    
   &   \hat{y} & = y  \\
    \hat{u} & =  \big( \underline{A_t}-C_t (\underline{A_t}^T) ^{-1} B_t \big)u + C_t (\underline{A_t}^T) ^{-1} v
   & \hat{v} & = -(\underline{A_t}^T) ^{-1} B_t u + (\underline{A_t}^T) ^{-1} v
\end{align*}
Let 
\begin{displaymath}
 \Omega = 
 \left (
 \begin{array}{cc}
     \underline{A_t}-C_t (\underline{A_t}^T) ^{-1} B_t   & \  C_t (\underline{A_t}^T) ^{-1}  \\
     -(\underline{A_t}^T) ^{-1} B_t  &  (\underline{A_t}^T) ^{-1}
 \end{array}
 \right )_{2m \times 2m}
\end{displaymath}
For $\xi \in V$, it admits a m-dimensional invariant torus $\mathbb{T}_{\xi} = \mathbb{T}^n \times \{0 \} \times \{ \xi \} \times \{0 \} $, $t\omega(\xi)$ is the rotational frequency.

In this paper, we consider the small perturbation, i.e.$P \neq 0$.

Set 
\begin{equation}
     \underline{A_t} = \sec{\theta^t} \  \    \   
    B_t = C_t  = \tan {\theta^t}
\end{equation}
where $\theta^t = (\theta_1^t, \theta_2^t, \dots, \theta_m^t)$ and
$\sec{\theta^t} = d i a g (\sec{\theta_1^t}, \sec{\theta_2^t}, \dots, \sec{\theta_m^t})$.
Then,
\begin{displaymath}
 \Omega = 
 \left (
 \begin{array}{cc}
     \cos{\theta^t}   & \   \sin{\theta^t}  \\
     - \sin{\theta^t}  & \   \cos{\theta^t}
 \end{array}
 \right )_{2m \times 2m}
\end{displaymath}
It is easy to see that $\Omega$ has the eigenvalues $e^{\pm i \theta_1^t}, e^{\pm i \theta_2^t} ,  \dots , e^{\pm i \theta_m^t}$.
Let us assume that the following non-resonant conditions hold for $l, i, j \in \mathbb{Z}$ and $1 \leq i, j \leq m$.
\begin{equation} \label{e04}
    \begin{cases}
         \theta_j^t - 2 \pi l \neq 0 \\
          \theta_j^t \pm  \theta_i^t - 2 \pi l \neq 0 , \  |l|+|i-j| \neq 0
    \end{cases}
\end{equation}

In this paper, we consider $h(\xi)$ is analytic in $\xi \in V + r$, where $V + r = \bigcup \limits_{b \in V} \{ \xi \in \mathbb{C}^n : \lvert \xi - b \rvert_2 < r \} \subseteq \mathbb{C}^n$, and 
$\omega(\xi) = h^{'}(\xi)$ satisfies the R\"{u}ssmann's non-degeneracy condition and satisfies 
\begin{equation}  \label{001}
    \lvert \omega(\xi_1) - \omega(\xi_2) \rvert \leq \Theta \lvert \xi_1 - \xi_2 \rvert , \  \xi_1 , \xi_2 \in V + r.
\end{equation}

\begin{remark} \label{re03}
For $\omega$ satisfies the R\"{u}ssmann's non-degeneracy condition, that is $\langle \omega , y \rangle \not\equiv 0 $ for all $y \in \mathbb{R}^n \setminus \{0\}$. 
In analytic cases, 
\cite{r5} gave an expression of the R\"{u}ssmann's non-degeneracy condition (see Remark 3.1 in \cite{r5}), that is, if the R\"{u}ssmann's non-degeneracy is satisfied, then $\exists \ \bar{n} \in \mathbb{N}$ such that 
 \begin{equation} 
     rank \big \{\partial^i_{\xi}\omega(\xi) :  |i| \leq \bar{n} \big \} =n ,  \   \forall \  \xi \in V
 \end{equation}
And, by R\"{u}ssmann (see Lemma 18.2 in \cite{r2}), there exist $ \bar{n} (V) \in \mathbb{N}, \beta_0 (V) >0 $ such that 
 \begin{equation}
     \min \limits_{\xi \in V} \max_{0 \leq v \leq \bar{n} } \lvert D^v \langle k, \omega(\xi) \rangle \rvert \geq \beta_0
 \end{equation}
where $\langle k,  \omega \rangle = \sum  \limits_{j=1}^n k_j \omega_j$, $\lvert k_j \rvert =1, j= 1, \dots, n $. And we can take the smallest of such integer  $\bar{n} $. 
And by R\"{u}ssmann (Theorem 18.4 in \cite{r2}), we have
\begin{equation}  
    \min_{\xi \in V} \max_{0\leq v\leq \bar{n}} \big \lvert D^v {\lvert k \rvert_2}^{-2}  | \langle k,  \omega(\xi) \rangle  |^2 \big \rvert \geq \beta_0
\end{equation}
for $\forall k \in \mathbb{Z}^n \setminus \{0 \}$. 
So, by $\lvert \cdot \rvert$ and $\lvert \cdot \rvert_2$ are equivalent, where $\lvert k \rvert = \sum  \limits_{j=1}^n \lvert k_j \rvert$, $\lvert k \rvert_2 = (\sum  \limits_{j=1}^n \lvert k_j \rvert^2)^{\frac{1}{2}} $, there exist $\bar{n}$, $\beta$ such that
\begin{equation}  \label{4}
    \min_{\xi \in V} \max_{0\leq v\leq \bar{n}} |D^v  \langle k,  \omega(\xi) \rangle |\geq \beta \lvert k \rvert
\end{equation}
for $\forall k \in \mathbb{Z}^n \setminus \{0 \}$, here $\bar{n} = \bar{n}(\omega, V) \in \mathbb{N} $  and $\beta = \beta(\omega, V) > 0$ are called the index and amount of $\omega$ with respect to $V$ respectively.
\end{remark}

Now, we give some notations. Set
\begin{equation*}
    T_s = \{ x \in \mathbb{C}^n / 2 \pi \mathbb{Z}^n  : | \Im x|_\infty \leq s \}
\end{equation*}
\begin{equation*}
    V_r = \{ y \in  \mathbb{C}^n : |y|_1 \leq r^2 \}   \  \  
    W_r = \{ w \in  \mathbb{C}^m : |w|_2 \leq r \}
\end{equation*}
Consider $D (s,r) = T_s \times W_r \times V_r \times W_r $. Here, $|x|_ \infty = \max  \limits_{1 \leq j \leq n} |x_j |, |y|_1 = \sum  \limits_ {j=1}^n |y_j|$, and $|w|_2 = ( \sum  \limits_ {j=1}^m |w_j|^2)^{\frac{1}{2}} $. 

Expand $P(x; \xi)$ into Fourier series 
$ P(x; \xi) = \sum \limits_{k \in   \mathbb{Z}^n} P_k (\xi)  e^{i \langle k, x \rangle} $, supposing  $P(x,u, \hat{y}, \hat{v}; \xi)$ is analytic on $D(s,r) $ with respect to $(x,u, \hat{y}, \hat{v} )$.
Define
\begin{equation*}
    \lVert P \rVert ^* _ s = \sum \limits_{k \in \mathbb{Z}^n } \lVert P_k \rVert ^* _V e^{s|k|}, with   \ 
    \lVert P_k \rVert ^* _V = \max \limits_{|i| \leq \bar{n} } \max \limits_ {\xi \in V} | \partial ^i _ \xi P_k(\xi) |
\end{equation*}
\begin{equation*}
    \lVert P \rVert ^* _ {s ,r} = \sum \limits_{k \in \mathbb{Z}^n } | P_k | ^* _r e^{s|k|} , with \ 
    | P_k | ^* _r = \sup \limits_{(u, \hat{y}, \hat{v}) \in W_r \times V_r \times W_r} \sum \limits_ {i,j,l} \lVert P_{k l i j} \rVert ^* _s {\hat{y}}^l u^i {\hat{v}}^j 
\end{equation*}
Define $X_p = (-\partial _{\hat{y}} P , -\partial _{\hat{v}} P, \partial_x P , \partial_u P)$. Define the weighed norm by 
\begin{equation*}
    |  \lVert  X_P  \rVert | _{r; D(s,r)} ^* = \lVert \partial _{\hat{y}} P \rVert _{s,r} ^* + \frac{1}{r} \lVert \partial _{\hat{v}} P \rVert _{s,r} ^* + \frac{1}{r^2} \lVert \partial _x P \rVert _{s,r} ^* + \frac{1}{r} \lVert \partial _u P \rVert _{s,r} ^*
\end{equation*}
Where $\lVert \partial _{\hat{x}} P \rVert _{s,r} ^* = \sum \limits _j \lVert \partial _{\hat{x}_j} P \rVert _{s,r} ^* $, $ \lVert \partial _{\hat{y}} P \rVert _{s,r} ^* = \max \limits_j \lVert \partial _{\hat{y}_j} P \rVert _{s,r} ^* $, and 
$ \lVert \partial _u P \rVert _{s,r} ^* = ( \sum \limits _j \lVert \partial _{u_j} P \rVert _{s,r} ^*)^ {\frac{1}{2}} $. $\lVert \partial _{\hat{v}} P \rVert _{s,r} ^* $ is similar to $\lVert \partial _u P \rVert _{s,r} ^* $.

Fix constant $\tau$ such that $\tau \geq (n+2)(\bar{n}+1)$ and $\tau \geq (n+2)L +1$, where $L$ is a constant in \ref{me1} Measure Estimation.

\section {Main Conclusions}
\label{se3}

For the symplectic mapping system with a small twist (\ref{e01}), we have two main theorems.
\begin{theorem} \label{t01}
    For the symplectic mapping $F$ defined by (\ref{e01}) , $t H = t N + t P$, assume $h$ and $P$ are real analytic on a complex neighborhood $D$ of  $ \mathbb{T}^n \times W \times V \times W $, $\omega $ satisfies the R\"{u}ssmann's non-degeneracy condition and (\ref{001}), and (\ref{e04}) holds, $\sup_{D} |P(x,u,\hat{y},\hat{v})| \leq \epsilon $ for $\epsilon > 0 $ small enough, and time step $t$ is small enough. Then, there is a Cantor-like subset $V_{\gamma} \subseteq V$ non-empty such that $F$ maintains a cluster of invariant torus $\{ \mathbb{T}_{\xi} \}_{\xi \in V_{\gamma}}$. Moreover, $meas(V \setminus V_{\gamma}) \to 0 $ as $\gamma \to 0 $.
\end{theorem}

\begin{theorem} \label{t02}
For the symplectic mapping $F$ defined by (\ref{e01}), and for the previous $t N, A_t, B_t, C_t$, suppose $\omega$  satisfies R\"{u}ssmann's non-degeneracy condition and (\ref{001}), $A_t, B_t, C_t$ satisfy the non-resonant condition (\ref{e04}), $P$ is real analytic on $D(s,r) \times V$ with $\xi \in V$, $t$ is small enough. Then, for $\gamma >0$,
$\exists \ \epsilon_0 >0$ small enough , such that, if $ |  \lVert  X_P  \rVert | _{r; D(s,r)} ^* = \epsilon \leq \epsilon_0$, then, $\exists \  V_{\gamma} \subseteq V$ Cantor-like and non-empty,  and there is a symplectic mapping $\Psi$, such that $F_ \infty = {\Psi}^{-1} \circ F \circ \Psi$ which is generated by $ t H_* = t N_* + t P_*$, where
\begin{align*} 
   t N_*(x,u,\hat{y},\hat{v};\xi) &= \langle t \omega_ \infty , \hat{y} \rangle +  \langle A_t^\infty  u , \hat{v} \rangle + \frac{1}{2}  \langle B_t^ \infty u , u \rangle + \frac{1}{2}  \langle C_t^ \infty \hat{v} , \hat{v} \rangle \\
     t P_*(x,u,\hat{y},\hat{v};\xi) &= t \sum \limits_{|i|+|j|+2|l| \geq 3} P_{l i j}(x;\xi) {\hat{y}}^l u^i {\hat{v}}^j
\end{align*}
with $\omega_ \infty - \omega = O(\epsilon)$, $A_t^ \infty - A_t = O(\epsilon)$, $B_t^ \infty - B_t= O(\epsilon)$, $\Psi - id = O(\epsilon)$. Then for each $\xi \in V_{\gamma}$, there exists an invariant torus $ \mathbb{T}_{\xi} $. Moreover, $meas(V \setminus V_{\gamma}) \to 0 $ as $\gamma \to 0$.
To be specific, 
   \begin{equation*}
     \lvert V \setminus V_{\gamma} \rvert \leq c_3   \gamma  d^{n}((n+1)^{-\frac{1}{2}}+2d+ \theta^{-1}d) + c_5   \gamma  d^{n} 
   \end{equation*}
where $c_3 = 12 \bar{c} (2\pi e)^\frac{n+1}{2} (\bar n +2)^{\bar n +3} [(\bar n  +2)!]^{-1} \big(\lvert \omega \rvert_{\mathcal{K}} + c  \epsilon +1 \big) \beta^{-\frac{\bar{n}+2}{\bar{n } +1}} \lvert  \omega \rvert^{\bar{ n}  +2}_{\mathcal{B}} $, $c_5 = 4  \bar{\bar{c}}  c_4 m (\lvert \omega \rvert_{V} + c  \epsilon + \max \limits_{1 \leq j \leq m}3|B^0_j |  +1) $, $\mathcal{B} = (V \times (0,1) + r) \cap \mathbb{R}^{n+1}$, $d$ is the diameter of $V \times (0,1)$, 
$\bar{c}$ and $\bar{\bar{c}}$ depend only on $n$, $\bar{n}$ and $\tau$,   $c_4$ and  $c$  are independent of $t$, $\epsilon$ and $\gamma$. 
\end{theorem}

Applying these conclusion to the symplectic difference scheme (\ref{004}), we get two main corollaries.
\begin{corollary} \label{co1}
 For the nearly integrable Hamiltonian system   (\ref{003}) and its corresponding symplectic difference scheme (\ref{004}) , $H_{\epsilon } = H_0 + \epsilon H_1$, $t H_0 = t N^{'} $ with the form (\ref{002}), assume $h$ is real analytic on a complex neighborhood $D$ of  $ \mathbb{T}^n \times W \times V \times W $, $\omega $ satisfies the R\"{u}ssmann's non-degeneracy condition and (\ref{001}), $A_t, B_t, C_t$ satisfy the non-resonant condition (\ref{e04}), $P$ is real analytic on $D(s,r) \times V$ with $\xi \in V$, $\sup_{D} |\epsilon H_1 + t^s P_1| \leq \epsilon^{'} $ for $\epsilon^{'} > 0 $. 
We can denote $\epsilon H_1 + t^s P_1 = P^{'}$. Then, for $\epsilon$ and time step $t$  small enough, most invariant tori of the nearly integrable system do not vanish, but are only slightly deformed, such that the symplectic difference scheme also has invariant tori.
 There exist a Cantor subset $V_{\epsilon,t} \subseteq V$, a mapping $\omega_{\epsilon,t} : V_{\epsilon,t} \to \Omega_{\epsilon,t}$, 
and a  symplectic mapping
$\Psi _{\epsilon,t} : V_{\epsilon,t} \times \mathbb{T}^n \to \mathbb{R}^n \times \mathbb{T}^n $, in the sense of Whitney, such that
 \begin{enumerate} 
    \item[(\romannumeral1)] 
    $G_{\infty} = \Psi _{\epsilon,t}^{-1} \circ G_{H_ \epsilon}^t \circ \Psi _{\epsilon,t} $ which is  generated by $ t (H_{\epsilon})_* = t N^{'}_* + t P^{'}_*$, where
\begin{align*} 
   t N^{'}_*(x,u,\hat{y},\hat{v};\xi) &= \langle t \omega_{\epsilon,t} , \hat{y} \rangle +  \langle A_t^\infty  u , \hat{v} \rangle + \frac{1}{2}  \langle B_t^ \infty u , u \rangle + \frac{1}{2}  \langle C_t^ \infty \hat{v} , \hat{v} \rangle \\
     t P^{'}_*(x,u,\hat{y},\hat{v};\xi) &= t \sum \limits_{|i|+|j|+2|l| \geq 3} P^{'}_{l i j}(x;\xi) {\hat{y}}^l u^i {\hat{v}}^j
\end{align*}
with $\omega_{\epsilon,t} - \omega = O(\epsilon^{'})$, $A_t^ \infty - A_t = O(\epsilon^{'})$, $B_t^ \infty - B_t= O(\epsilon^{'})$,$\Psi_{\epsilon,t} - id = O(\epsilon^{'})$. 
Then for each $\xi \in V_{\epsilon,t}$, there exists an invariant torus $ \mathbb{T}_{\xi} $.
    \item[(\romannumeral2)] 
    $meas(V \setminus V_{\epsilon,t}) \to 0 $ as $\gamma \to 0$.
    To be specific, 
    $V_{\epsilon,t}$ is a set of positive measure with 
    the measure estimate
    \begin{equation}
        \lvert V \setminus V_{\epsilon,t} \rvert \leq c_3^{'}   \gamma  d^{n}((n+1)^{-\frac{1}{2}}+2d+ \theta^{-1}d) + c_5^{'}  \gamma  d^{n}
    \end{equation}
    where $c_3^{'} = 12 \bar{c} (2\pi e)^\frac{n+1}{2} (\bar n +2)^{\bar n +3} [(\bar n  +2)!]^{-1} \big(\lvert \omega \rvert_{\mathcal{K}} + c  \epsilon^{'} +1 \big) \beta^{-\frac{\bar{n}+2}{\bar{n } +1}} \lvert  \omega \rvert^{\bar{ n}  +2}_{\mathcal{B}} $, $c_5^{'} = 4  \bar{\bar{c}}  c_4 m (\lvert \omega \rvert_{V} + c  \epsilon^{'} + \max \limits_{1 \leq j \leq m}3|B^0_j |  +1) $,
    $\mathcal{B}$, $d$, 
    $\bar{c}$, $\bar{\bar{c}}$ , $c_4$ and  $c$  are in Theorem \ref{t02}. 
\end{enumerate}
\end{corollary}
\begin{corollary} \label{co2}

         Let the time step be $t_1$ and $t_2$ respectively, and comppare the two systems, then for $\xi \in V_{\epsilon,t_1} \cap V_{\epsilon,t_2}$ we have

\begin{align*}
    & \lvert \Psi _{\epsilon,t_1} - \Psi _{\epsilon,t_2} \rvert \leq c (t_1^{\frac{s}{4}} - t_2^{\frac{s}{4}}), \ \lvert \omega_{\epsilon,t_1} - \omega_{\epsilon,t_2} \rvert \leq c (t_1^{\frac{s}{2}} - t_2^{\frac{s}{2}}) .
\end{align*}
    $ V_{\gamma,t_1} \cap V_{\gamma,t_2} $ is also a set of positive measure if $\gamma$ is small enough. And for $\gamma \to 0$, we also have $\lvert V \setminus (V_{\gamma,t_1} \cap V_{\gamma,t_2}) \rvert \to 0 $.
\end{corollary}

\section {The Symplectic Mapping System with a Small Twist (\ref{e01})}
\label{se4}

Let us focus on the symplectic mapping system with a small twist (\ref{e01}). 
The proof of Theorem \ref{t01}  and Theorem \ref{t02} are similar to the proof of Theorem 1.1 and Theorem 1.2 in \cite{r6}, so we will not write them in detail here. we are mainly concerned with measure estimation.

The symplectic mapping in  \cite{r6} is expressed as 
\begin{align*} 
      \hat{x} & = x+  \partial_{\hat{y}}H(x,u,\hat{y},\hat{v})     & 
    y & = \hat{y} +  \partial_{x}H(x,u,\hat{y},\hat{v})  \\
    \hat{u} & = u+ \partial_{\hat{v}}H(x,u,\hat{y},\hat{v})      & 
    v & = \hat{v} + \partial_{u}H(x,u,\hat{y},\hat{v})
\end{align*}
where $H= N+P$, and 
\begin{equation*}
     N(x,u,\hat{y},\hat{v}) =   h(\hat{y}) + \langle (A-I) u,\hat{v} \rangle +  \frac{1}{2} \langle B u,u \rangle + \frac{1}{2} \langle C \hat{v},\hat{v} \rangle
\end{equation*}
We introduced a positive step $t$, and we want to explore the system as $t \to 0$.

Now, we briefly describe the proving process of the Theorem \ref {t01} and  Theorem \ref {t02}, please refer to \cite{r6} for details.
Write $x_+$ as $x_{v+1}$ for simplicity, and everything else is similar.

Let $p = (x,u), q = (y,v)$. We have $F:(p,q) \to (\hat{p},\hat{q})$
\begin{equation}
       \hat{p}  = p + \partial_{\hat{q}} t  H(p,\hat{q}) \  \    \   
        q  = \hat{q} + \partial_{p} t  H(p,\hat{q})
\end{equation}
Suppose that $\Psi$ is close to the identity, $\Psi : (p_+,q_+) \to (\hat{p},\hat{q}) $
\begin{equation} \label{e07}
       q  = q_+ + F_1(p, q_+) \  \    \   
       p_+ = p + F_2(p, q_+)
\end{equation}
then we have $F_+ = {\Psi}^{-1} \circ F \circ \Psi : (p_+, q_+) \to (\hat{p}_+, \hat{q}_+)$, and it can be expressed implicitly by 
\begin{equation}
       \hat{p}_+  = p + t H_2(p, \hat{q}) + F_2(\hat{p}, \hat{q}_+) \  \    \   
       q_+ = \hat{q} +t H_1(p, \hat{q} ) - F_1(p, q_+)
\end{equation}

Then, by a similar but slightly different processing from Lemma 2.1 in \cite{r6}, $F_+$ can be expressed implicitly as 
\begin{equation}
       \hat{p}_+  = p_+ + \partial_{\hat{q}_+} t  H_+(p_+,\hat{q}_+) \  \    \   
        q_+  = \hat{q}_+ + \partial_{p_+} t  H_+(p_+,\hat{q}_+)
\end{equation}
where $t H_+$ has the following form
\begin{align*}
    t H_+(p_+,\hat{q}_+) = & \ t H(p, \hat{q}) + t H_1(p, \hat{q}) F_2(p, q_+) - t H_2(p, \hat{q}) F_1(\hat{p}, \hat{q}_+) 
    + F(\hat{p}, \hat{q}_+) \\
    & - F(p, q_+) - F_1(p, q_+) F_2(p, q_+) + F_1(\hat{p}, \hat{q}_+) F_2(p, q_+)
\end{align*}
And for $t H = t N + t P$, by a similar analysis to  \cite{r6}, let $z = (p_+, \hat{q}_+)$ then  we can get 
\begin{equation} \label{e010}
    t H_+(z) = t H(z) + F(p_+ + t N_2(z), \hat{q}_+ ) - F(p_+, \hat{q}_+ + t N_1 (z) ) + Q (z)
\end{equation}
where $Q$ is a higher order small quantity.

\begin{remark}
To be specific, by some similar but slightly different analysis from  \cite{r6}, we have these following expressions. 
From the foregoing we have let $z = (p_+, \hat{q}_+)$, here we let $\Delta z = (-F_2(p,q_+), F_1(\hat{p}, \hat{q}_+))$, then $(p, \hat{q}) = z+ \Delta z$. So,
\begin{equation*}
    t H(p, \hat{q}) = t H(z+ \Delta z) = t H(z) + \langle \nabla t H (z), \Delta z \rangle - Q_1
\end{equation*}
where 
\begin{equation*}
    Q_1 = -\frac{1}{2} \int_0^1 \langle D^2 t H (z_s) \Delta z, \Delta z \rangle d s, \   z_s = z + s \Delta z
\end{equation*}
here $\nabla t H$ means the gradient of $t H$, and $D^2 t H$ indicates the Hessian matrix. By
\begin{equation*}
    \nabla t H (z+ \Delta z) = \nabla t H(z) + \int_0^1 D^2 t H(z_s) \Delta z d s
\end{equation*}
we have
\begin{align*}
    & t H_1(p, \hat{q}) F_2(p, q_+) - t H_2(p, \hat{q}) F_1(\hat{p}, \hat{q}_+) \\
    = & - \langle \nabla t H(p, \hat{q}), (-F_2(p,q_+), F_1(\hat{p}, \hat{q}_+)) \rangle \\
    = & - \langle \nabla t H (z+ \Delta z), \Delta z \rangle 
    = - \langle \nabla t H (z), \Delta z \rangle +2 Q_1
\end{align*}
Note that 
\begin{align*}
    \hat{p} = p_+ + t N_2(z + \Delta z) +t P_2(p, \hat{q}) - F_2(p,q_+)\\
    t N_2(z + \Delta z) = t N_2(z) + \int_0^1 D t N_2(z_s) \Delta z d s
\end{align*}
then we have
\begin{align*}
    &  F(p_+ + t N_2(z), \hat{q}_+) - F(\hat{p}, \hat{q}_+)  \\
    = & \  Q_2 + 
     F_1(\hat{p} , \hat{q}_+)  F_2(p, q_+)+ \sum \limits_{i=2}^\infty \frac{1}{i !}  F_{1^i}(\hat{p}, \hat{q}_+) \cdot \big( F_2(p, q_+)  \big) ^i
\end{align*}
where 
\begin{equation*}
    Q_2 = \int_0^1 F_1(z+ s \Delta w ) \Delta w d s  , \ 
    \Delta w = t P_2(p, \hat{q}) + \int_0^1 D t N_2(z_s) \Delta z d s 
\end{equation*}
where $D t N_2$ indicates the Jacobian matrix, and $F_{1^i}(\hat{p}, \hat{q}_+)$ means the $i$-th partial of $F(\hat{p}, \hat{q}_+)$ with respect to $\hat{p}$. Similarly,  by $p_+ = p + F_2(p, q_+)$, we have
\begin{align*}
    F(p_+, q_+) - F(p, q_+) & = F_1(p, q_+) F_2(p, q_+) + \sum \limits_{i=2}^\infty \frac{1}{i !}  F_{1^i} (p, q_+) \big( F_2(p, q_+) \big)^i \\
    & F(p_+, \hat{q}_+ + t N_1(z)) - F(p_+, q_+) = Q_3 
\end{align*}
where
\begin{equation*}
    Q_3 = \int_0^1 F_2(z+ s \Delta \omega^{'}) \Delta \omega^{'} d s , \ \Delta \omega^{'} = \partial_1 \big(t H_+(z) - t H(z) + t P(z) \big) 
\end{equation*}
And we have
\begin{equation*}
    \sum \limits_{i=2}^\infty \frac{1}{i !}  F_{1^i}(\hat{p}, \hat{q}_+) \cdot \big( F_2(p, q_+)  \big) ^i 
    - \sum \limits_{i=2}^\infty \frac{1}{i !}  F_{1^i} (p, q_+) \big( F_2(p, q_+) \big)^i = Q_4
\end{equation*}
where
\begin{equation*}
    Q_4 = \big( F (\hat{p}, \hat{q}_+) - F (p, q_+) \big) - \big( F_1 (\hat{p}, \hat{q}_+) - F_1 (p, q_+)  \big) F_2(p, q_+)
\end{equation*}
So, $Q = Q_1 - Q_2 + Q_3 - Q_4$.
\end{remark}

Next, similar to \cite{r6}, we truncate $P$ with $P = R +  \widetilde{P} $, where 
\begin{align}  \label{e011}
    R(p, \hat{q}) = P_{000}(x) + \langle P_{100}(x), \hat{y} \rangle + \langle P_{010}(x), u \rangle + \langle P_{001}(x), \hat{v} \rangle  \notag \\
    + \langle P_{011}(x)u, \hat{v} \rangle + \frac{1}{2} \langle P_{020}(x)u, u \rangle + \frac{1}{2} \langle P_{002}(x) \hat{v}, \hat{v} \rangle 
\end{align}
with the coefficients
$ P_{l i j} = \frac{\partial ^{l+i+j} P }{\partial{{\hat{y}^l} \partial{u^i} \partial{{\hat{v}^j} }}} \big | _{u=0, \hat{y}=0, \hat{v}=0} $, 
for $ 2l + i + j \leq 2 $.
Then (\ref{e010}) leads to the following equation:
\begin{align*}
      t R(p_+, \hat{q}_+)  - & F(p_+, \hat{q}_+ + t N_p (p_+, \hat{q}_+) )  + F(p_+ + t N_q(p_+, \hat{q}_+), \hat{q}_+) \\
        & + t N(p_+, \hat{q}_+) 
        = t \bar{N}(p_+, \hat{q}_+) 
\end{align*}
Where $ t \bar{N}(p_+, \hat{q}_+) $ represent the new normal form.

Then, similar to \cite{r6}, let us use a linear operator, which is defined as 
$L F = F(p + t N_{\hat{q}}(p, \hat{q}), \hat{q}) - F(p, \hat{q} + t N_p (p, \hat{q}) ) $. 
Let $F(p, \hat{q})$ remain in the form of $R(p, \hat{q})$ in (\ref{e011}), then we have
\begin{align*}
     F(p, \hat{q} + t N_p (p, \hat{q}) ) & =  F_{000}(x) + \langle F_{100}(x), \hat{y} \rangle + \langle F_{010}(x), u \rangle + \langle B_t F_{001} (x), u \rangle \\
     & + \langle  \underline{A_t} F_{001} (x), \hat{v}   \rangle + \langle  \underline{A_t} F_{011} (x) u , \hat{v}   \rangle + \langle  \underline{A_t} F_{002} (x) B_t u, \hat{v}   \rangle \\
     & +  \langle B_t F_{011} (x) u, u \rangle + \frac{1}{2} \langle B_t F_{002} (x) B_t u, u \rangle \\
     & + \frac{1}{2} \langle \underline{A_t} F_{002} (x) {\underline{A_t}}^T \hat{v}, \hat{v} \rangle + \frac{1}{2}  \langle F_{020} (x) u, u \rangle
\end{align*}
\begin{align*}
     F(p + t N_{\hat{q}}(p, \hat{q}), \hat{q}) & =  F_{000}(\tilde{x}) + \langle F_{100}(\tilde{x}), \hat{y} \rangle + \langle F_{010}(\tilde{x}), \hat{v} \rangle + \langle C_t F_{001} (\tilde{x}), \hat{v} \rangle \\
     & + \langle  {\underline{A_t}}^T F_{010} (\tilde{x}), u  \rangle + \langle  F_{011} (\tilde{x})  \underline{A_t} u , \hat{v}   \rangle + \langle  C_t F_{020} (\tilde{x}) \underline{A_t} u, \hat{v}   \rangle \\
     & +  \langle F_{011} (\tilde{x}) C_t  \hat{v}, \hat{v} \rangle + \frac{1}{2} \langle C_t F_{020} (\tilde{x}) C_t \hat{v}, \hat{v} \rangle \\
     & + \frac{1}{2} \langle  {\underline{A_t}}^T F_{020} (\tilde{x}) \underline{A_t} u, u  \rangle + \frac{1}{2}  \langle F_{002} (\tilde{x}) \hat{v}, \hat{v} \rangle
\end{align*}
where $ \tilde{x} = x + t \omega $ , $R_{020}, R_{002} , F_{020}$ and $F_{002}$ are symmetric. So
\begin{equation*}
    L F = L_0 + \langle L_1, \hat{y} \rangle + \langle L_{10}, u \rangle + \langle L_{01}, \hat{v} \rangle +  \langle L_{11} u , \hat{v} \rangle + \frac{1}{2}  \langle L_{20} u, u \rangle +  \frac{1}{2}  \langle L_{02} \hat{v}, \hat{v} \rangle
\end{equation*}
where 
\begin{equation} \label{e041}
      L_0 = F_{000} (\tilde{x}) - F_{000} (x)\  \    \  
      L_1 = F_{100} (\tilde{x}) - F_{100} (x) 
\end{equation}
\begin{equation} \label{e042}
    \begin{cases}
        L_{10} = {\underline{A_t}}^T F_{010} (\tilde{x}) - F_{010} (x) - B_t F_{001} (x) \\
        L_{01} = C_t F_{010} (\tilde{x}) + F_{001} (\tilde{x}) - \underline{A_t} F_{001} (x)
    \end{cases}
\end{equation}
\begin{equation}  \label{e043}
    \begin{cases}
        L_{11} = F_{011} (\tilde{x}) \underline{A_t}  - \underline{A_t} F_{011} (x) + C_t F_{020} (\tilde{x}) \underline{A_t} - \underline{A_t} F_{002} (x) B_t \\
        L_{11}^T = {\underline{A_t}}^T F_{011}^T (\tilde{x}) -  F_{011}^T (x) {\underline{A_t}}^T + {\underline{A_t}}^T F_{020} (\tilde{x}) C_t - B_t F_{002} (x) {\underline{A_t}}^T \\
        L_{20} = {\underline{A_t}}^T F_{020} (\tilde{x}) \underline{A_t} - F_{020} (x) - B_t F_{002} (x) B_t - B_t F_{011} (x) - F_{011}^T (x) B_t \\
        L_{02} = C_t  F_{020} (\tilde{x}) C_t +  F_{020} (\tilde{x}) - \underline{A_t} F_{002} (x) {\underline{A_t}}^T + F_{011} (\tilde{x}) C_t + C_t  F_{011}^T (\tilde{x})
    \end{cases}
\end{equation}
Then we solve the following equations
\begin{equation}  \label{e012}
    \begin{cases}
        L_0 = t R_{000} (x) - [ t R_{000} ], 
        L_1 = t R_{100} (x) - [ t R_{100} ]  \\
        L_{10} = t R_{010} (x),    
        L_{01} = t R_{001} (x) \\
        L_{11} = t R_{011} (x) - t \hat{A},
        L_{11}^T = t R_{011}^T (x) - t {\hat{A}}^T,  \\
        L_{20} = t R_{020} (x) - t \hat{B},
        L_{02} = t R_{002} (x) - t \hat{C}
    \end{cases}
\end{equation}
where $\hat{A} , \hat{B} $ and $\hat{C}$ will be explained later, $[\cdot]$ represents the mean value.
In order to avoid irreversible of the coefficient matrix of the system, we need to satisfy the small divisor conditions listed as:
\begin{equation}\label{e031}
    |\langle k,t\omega_v \rangle -2\pi l| \geq \frac{t\gamma_v}{({1+|k|)}^{\tau}} 
\end{equation} 
\begin{equation} \label{e032}
    |\langle k,t\omega_v \rangle -\theta^{t, v}_j -2\pi l| \geq \frac{t\gamma_v}{({1+|k|)}^{\tau}} 
\end{equation}
\begin{equation} \label{e033}
    |\langle k,t\omega_v \rangle -\theta^{t, v}_i -\theta^{t, v}_j -2\pi l| \geq \frac{t\gamma_v}{({1+|k|)}^{\tau}} 
\end{equation}
\begin{equation} \label{e034}
    |\langle k,t\omega_v \rangle +\theta^{t, v}_i -\theta^{t, v}_j -2\pi l| \geq \frac{t\gamma_v}{({1+|k|)}^{\tau}} 
\end{equation}
 for all $ \ k\in \mathbb{Z}^n \setminus \{ 0\} $, $v= 0,1,2,3, \dots$ and $i \neq j$, where $\tau$ is a fixed integer.
 As for $k = 0$, the corresponding small divisor conditions should also be satisfied, as detailed in Remark \ref{re11} of the measure estimation section. 
Then we can solve for the coefficients of $F$, and we get $F$.
\begin{remark} \label{re003}
Similar to \cite{r6}, as $\underline{A_t}$, $B_t$, $C_t$ are all diagonal, the operator $L:(F_{010}, F_{001}) \to (L_{10}, L_{01})$ decomposes into $m$ independent operators. 
Now let $(M_k)_{2 \times 2}$ represent the coefficient matrix of (\ref{e042}), similar to \cite{r6} we have $| \det (M_k) | \geq \frac{c t^2 \gamma^2}{(|k|+1)^{2\tau}}$.
By the result: $ | \det A | \leq c(n) l^n$ , where $A_{n\times n} = (a_{i j})$ is a matrix , $ l = \max \limits_{1\leq i,j \leq n} |a_{i j}|$. Then,  $\sum \limits_{1\leq i,j \leq n} |(M_k)_{i j}| \geq \max \limits_{1\leq i,j \leq n} |(M_k)_{i j}| \geq \frac{c t \gamma}{(|k|+1)^{\tau}}$. So similar to \cite{r6}, $\lVert X \rVert_{s-\rho}^* \leq \frac{c \lVert Y \rVert_s^*}{\gamma^{\bar{n}+1} \rho^{\tau (\bar{n}+1)+\bar{n}+n}}$ with $X=(F_{010}, F_{001})^T$ and $Y=(R_{010}, R_{001})^T$.
The coefficient matrix of (\ref{e043}) is similar.
\end{remark}

Let 
\begin{align*}
    \hat{\omega} & = d i a g (\hat{\omega}_1 , \dots , \hat{\omega}_m) & \hat{A} & = d i a g (\hat{A}_1 , \dots , \hat{A}_m) \\
    \hat{B} & = d i a g (\hat{B}_1 , \dots , \hat{B}_m) & 
    \hat{C} & = d i a g (\hat{C}_1 , \dots , \hat{C}_m)
\end{align*}
with 
\begin{equation*}
    \hat{\omega}_j = [R_{100}^{j j}]  \  \  \hat{A}_j = [R_{011}^{j j}]  \  \    \hat{B}_j = [R_{020}^{j j}]   \  \    \hat{C}_j = [R_{002}^{j j}]
\end{equation*}
Let 
\begin{equation*}
   \hat{N}(p, \hat{q}) = \langle t \hat{\omega}, \hat{y} \rangle + \langle t \hat{A} u, \hat{v} \rangle  + \frac{1}{2} \langle t \hat{B} u, u \rangle  + \frac{1}{2} \langle t \hat{C} \hat{v}, \hat{v} \rangle 
\end{equation*}
Then equation (\ref{e012}) is $L F = t R - \hat{N}$, similar to \cite{r6} we have
\begin{equation*}
    | \lVert X_F   \rVert | _{r; D(s- \rho , r)} ^* \leq \frac{c   \epsilon}{ \gamma^{\bar{\nu}} \rho ^ \nu} 
\end{equation*}
where $\bar{\nu} = \bar{n}+1$ , $\nu =  \tau (\bar{n}+1)+ \bar{n}+n$.
 And $\psi$ can be expressed as
\begin{equation*}
    \psi (p, q) = (p_+, q_+) + X_F (p, q_+) \ with  \ X_F = (-F_{q_+}, F_p)
\end{equation*}
By the same calculation method  in \cite{r6} , we have $\psi : D(s-2 \rho , r/2) \to D(s- \rho , r)$, moreover 
\begin{equation*}
    |  \lVert  \psi - id  \rVert  | _{r; D(s-3 \rho , r/4)} ^* \leq \frac{c   \epsilon}{ \gamma^{\bar{\nu}} \rho ^ \nu} 
    \ 
    |  \lVert D \psi - id \rVert | _{r,r; D(s-3 \rho , r/4)} ^* \leq \frac{c   \epsilon}{ \gamma^{\bar{\nu}} \rho ^ {\nu +1}} 
\end{equation*}
where $| \lVert \cdot \rVert |_{r,r}$ represent the operator norm under $| \lVert \cdot \rVert | _r$.

Similar to \cite{r6} , now let $\bar{F} = {\psi}^{-1} \circ F \circ \psi $, then it can be generated by $t \bar{H} = t \bar{N} + t \bar{P}$ with $\bar{P} = \widetilde{P} + Q $, and 
\begin{equation*}
    t \bar{N} = t N + t \hat{N} = \langle t \omega_+, \hat{y} \rangle + \langle \bar{A_t} u,\hat{v} \rangle +  \frac{1}{2} \langle \bar{B_t} u,u \rangle + \frac{1}{2} \langle \bar{C_t} \hat{v},\hat{v} \rangle
\end{equation*}
where $t \omega_+ = t \omega + t \hat{\omega}$ , $\bar{A_t} = A_t + t \hat{A}$, $\bar{B_t} = B_t + t \hat{B}$, $\bar{C_t} = C_t + t \hat{C}$ with $\lVert \hat{\omega} \rVert , \lVert \hat{A} \rVert , \lVert \hat{B} \rVert , \lVert \hat{C} \rVert \leq \epsilon$.

Note that $\underline{\bar{A_t} }^2 - \bar{B_t} ^2 = I $ and $\bar{B_t} = \bar{C_t}$ may no longer hold with $\underline{\bar{A_t} } = I + \bar{A_t}$, so we need normalization. 
There exists $\bar{\psi}$ symplectic, and it can transform $t \bar{N}$ into a normal form
\begin{equation*}
    t N_+ = {\bar{\psi}}^{-1} \circ \bar{F} \circ \bar{\psi} =  \langle t \omega_+, \hat{y} \rangle + \langle A_t^+ u,\hat{v} \rangle +  \frac{1}{2} \langle B_t^+  u,u \rangle + \frac{1}{2} \langle C_t^+ \hat{v},\hat{v} \rangle
\end{equation*}
where $I + A_t^+ =  \underline{A_t^+} = \sec \theta ^t_+$, $B_t^+ = C_t^+ = \tan \theta ^t_+$, and $\lVert \theta ^t_+ - \theta^t \rVert ^* \leq c t  \epsilon $.
Let $\Psi = \psi \circ \bar{\psi}$, similar to \cite{r6}, $\Psi$ is generated by
\begin{align*}
    \langle x, y_+ \rangle +  \langle \lambda u , v_+ \rangle +  F_{000}(x) +  \langle F_{100}(x), y_+ \rangle + \langle F_{010}(x), u \rangle +\langle  \lambda F_{001}(x), v_+ \rangle \\
    +  \langle  \lambda F_{011}(x) u, v_+  \rangle +  \frac{1}{2} \langle F_{020}(x) u, u \rangle + \frac{1}{2} \langle \lambda \beta v_+ , v_+ \rangle +  \frac{1}{2} \langle {\lambda}^2 F_{002}(x) v_+ , v_+ \rangle
\end{align*}
Where $\lambda = \lambda (t)$ and $\beta = \beta (t)$, please refer to the normalization process in \cite{r6} (Lemma 3.5 in \cite{r6}), and $\lambda - I = O (\epsilon)$ , $\beta = O(\epsilon) $. We have
$t H_+ = t H \circ \Psi = t N_+ + t P_+$, and
\begin{equation*}
  | \lVert  \Psi - id \rVert | _{r; D(s-3 \rho , r/4)} ^* \leq \frac{c   \epsilon}{ \gamma^{\bar{\nu}} \rho ^ \nu}  
    \ 
    | \lVert D \Psi - id \rVert | _{r,r; D(s-3 \rho , r/4)} ^* \leq \frac{c   \epsilon}{ \gamma^{\bar{\nu}} \rho ^ {\nu +1}} 
\end{equation*}

As for $t P_+$, similar to \cite{r6}, let us focus on the  truncated remainder $t \widetilde{P}$ and the higher order small quantity $Q$. Then we have $ | \lVert X_{t \widetilde{P}} \rVert |  _{\eta r; D(s-5 \rho , \eta r)} ^* \leq \eta t \epsilon$, 
and $ | \lVert X_{Q} \rVert |  _{r; D(s-5 \rho , r/16)} ^* \leq  \frac{c t \epsilon^2}{ \gamma^{2 \bar{\nu}} \rho ^ {2 \nu }} $. As $t \bar{P} = t \widetilde{P} +Q$, 
we have $ | \lVert X_{t \bar{P}} \rVert |  _{\eta r; D(s-5 \rho , \eta r)} ^* \leq \eta t \epsilon +   \frac{c \eta^{-2} t \epsilon^2}{ \gamma^{2 \bar{\nu}} \rho ^ {2 \nu }} $. 
Let $\eta ^3 = \frac{\epsilon}{ \gamma^{\bar{\nu}} \rho ^ {\nu}}$, $r_+ = \eta r $, $s_+ = s - 5 \rho $  and $\epsilon _+ = c \eta \epsilon$,
then 
\begin{equation*}
    | \lVert X_{t P_+} \rVert |  _{r_+; D(s_+ , r_+)} ^* \leq t \epsilon_+
\end{equation*}

Consider KAM iteration, similar to \cite{r6}, let $\Psi _0 = \Psi, \omega_{t0} = \omega_t, \theta^{t,0} = \theta^t$, and $s_0= s,  \rho _0 = s/20 , r_0= r, \epsilon _0 = \epsilon, \gamma_0 = \gamma $, define
\begin{equation*}
    \rho _v = \frac{\rho _0}{2^v} , \gamma_v = \frac{\gamma^{(\bar{n}+1) L}}{2^{(\bar{n}+1)L v}} , \eta _v =  \frac{ \epsilon _v}{ \gamma_v^{2{\bar{\nu}}} \rho_v ^ {2\nu}}, r_{v+1} = \eta_v r_v, s_{v+1} = s_v -5 \rho_v
\end{equation*}
Where $L$ is a fixed integer, see latter Lemma \ref{l2}, then 
\begin{equation*}
    \lVert t \omega_{v+1} - t \omega_v \rVert ^* \leq c t  \epsilon_v ,  \lVert \theta^{t , v+1} - \theta^{t, v} \rVert ^* \leq c t \epsilon_v ,  | \lVert X_{ t P_v}   \rVert | _{r_v; D(s_v , r_v)} ^* \leq t \epsilon_v
\end{equation*}
for $\forall \ j \geq 1$. Let $F_{v+1} = ({\Psi}^v)^{-1} \circ F_v \circ {\Psi}^v$ with
 ${\Psi}^v = {\Psi}_0 \circ {\Psi}_1 \circ \dots \circ {\Psi}_v $, we have
 \begin{align*}
        &  | \lVert \Psi_v - id   \rVert | _{r_v; D(s_v-3 \rho_v , r_v)} ^* \leq \frac{c   \epsilon_v}{ \gamma_v^{\bar{\nu}} {\rho_v} ^ \nu}  \\
      & | \lVert D \Psi_v - id  \rVert | _{r_v,r_v; D(s_v-3 \rho_v , r_v)} ^* \leq \frac{c   \epsilon_v}{\gamma_v^{\bar{\nu}} {\rho_v} ^ {\nu +1}}  
 \end{align*}
Set $\Psi_{\infty} = \lim \limits _{v \to \infty} {\Psi}^v$, let $t H_v = t N_v + t P_v \to t H_* = t N_* + t P_*$ with
\begin{equation*}
    t N_* = \langle  t \omega_{\infty}(\xi), \hat{y} \rangle + \langle A_t ^\infty (\xi) u , \hat{v}  \rangle + \frac{1}{2} \langle B_t^ \infty (\xi) u , u  \rangle +  \frac{1}{2} \langle C_t ^\infty (\xi) \hat{v} , \hat{v}  \rangle
\end{equation*}
where $I + A_t ^\infty = \underline{A_t ^\infty} = \sec {\theta^{t, \infty}}$ , $B_t ^\infty = C_t  ^ \infty = \tan {\theta^{t, \infty}}$ . And
\begin{equation*}
    \lVert t \omega_{\infty} - t \omega \rVert ^* \leq c t  \epsilon ,  \lVert \theta^{t , \infty} - \theta^t \rVert ^* \leq c t \epsilon 
\end{equation*}

Thus, we get a family of invariant torus for $\xi \in V_\gamma$. 
As for measure estimation in Theorem \ref{t01}, next section gives it.

\begin{remark}  \label{re10}
By Lemma 3.5 in \cite{r6}, we have
\begin{equation*}
    \underline{A_t^+} = \frac{\lambda \underline{A_t} }{\lambda + \beta B_t}, 
    B_t^+ = C_t^+ = \frac{B_t}{\lambda (\lambda + \beta B_t)}
\end{equation*}
By $\lambda - I = O (\epsilon)$ , $\beta = O(\epsilon) $, for $\epsilon$ small enough  we have $\frac{|B_t^+|}{|B_t|} \leq \frac{1}{1- c \epsilon}$, where $c$ is a positive constant. 
And by $\epsilon_+ \sim \epsilon^{\frac{4}{3}}$, so for $\epsilon$  small enough, we have $c \epsilon_v \leq \frac{1}{2^{v+2}-1}$. 
Then 
$\frac{| B_{t}^{v+1}|}{| B_t^{v} |} \leq \frac{2^{v+2}-1}{2^{v+2}-2} = \frac{\sum \limits_{i=0}^{v+1} \frac{1}{2^i}}{\sum \limits_{i=0}^{v} \frac{1}{2^i}}$. 
By induction, if $| B_t^{v} | \leq \sum \limits_{i=0}^{v} \frac{1}{2^i}$, then $| B_t^{v+1} | \leq \sum \limits_{i=0}^{v+1} \frac{1}{2^i}$.
So we have  $\frac{| B_t^v |}{| B_t^0 |} \leq \sum \limits_{i=0}^{v} \frac{1}{2^{i}} \leq 2$, for $B_t^v =  B^v t $, we have $ | B^v t | \leq 2 | B^0 t |$ for $v = 0, 1,2,3 \dots$.
Similarly, we have $\frac{|B_t^+|}{|B_t|} \geq \frac{1}{1 + c \epsilon}$, and we can get $c \epsilon_v \leq 1$, then we have $\frac{|B_t^{v+1}|}{|B_t^v|} \geq \frac{1}{2} = \frac{2- \sum \limits_{i=0}^{v+1} \frac{1}{2^i}}{2 - \sum \limits_{i=0}^{v} \frac{1}{2^i}}$. So by induction, we have 
$| B_t^v | \geq 2 - \sum \limits_{i=0}^{v} \frac{1}{2^i}$, so
$\frac{| B_t^v |}{| B_t^0 |} \geq 2 - \sum \limits_{i=0}^{v} \frac{1}{2^{i}} \geq \frac{1}{2}$, $ | B^v t | \geq \frac{1}{2} | B^0 t |$ for $v = 1,2,3 \dots$.
This remark is in preparation for the later measure estimation.
\end{remark}

\section {Measure estimation} \label{me1}
Firstly, for all $ \ k\in \mathbb{Z}^n \setminus \{0\} $ and $i \neq j$, let us focus on the following four conditions.
\begin{equation}\label{e1}
    |\langle k,t\omega_v(\xi) \rangle -2\pi l| \geq t\gamma_v |k|^{-\tau}
\end{equation} 
\begin{equation} \label{e2}
    |\langle k,t\omega_v(\xi) \rangle -\theta^{t, v }_j -2\pi l| \geq t\gamma_v |k|^{-\tau}
\end{equation}
\begin{equation} \label{e3}
    |\langle k,t\omega_v(\xi) \rangle -\theta^{t, v }_i -\theta^{t, v }_j -2\pi l| \geq t\gamma_v |k|^{-\tau}
\end{equation}
\begin{equation} \label{e4}
    |\langle k,t\omega_v(\xi) \rangle +\theta^{t, v }_i -\theta^{t, v }_j -2\pi l| \geq t\gamma_v |k|^{-\tau}
\end{equation}
   We want  $t \omega_v $ satisfy them with corresponding $\gamma_v$ for $v=0,1,2,3 \dots$

Firstly, 
in order for $t \omega_v, v=0,1,2,3 \dots$ to satisfy (\ref{e1}), 
we focus on  
\begin{equation*}
     \big \{\xi \in V : \lvert  \langle k, t \omega_v(\xi) \rangle -2 \pi l  \rvert \geq \frac{t \gamma_v}{{\lvert k \rvert}^\tau} , \forall \  k \in \mathbb{Z}^n \setminus \{0 \}, \forall  \ l \in \mathbb{Z} \big \}
\end{equation*}
Then, we know that, if $\xi$ is a member of the set, then (\ref{e1}) holds for such $\xi$. So we estimate its measure, and we give the following theorem.

\begin{theorem} \label{t1}
    Let $ V \subseteq \mathbb{R}^n $ be compact , $d$ is the diameter of $V \times (0,1)$,
   define $\mathcal{B}:= {V \times (0,1) + r }\subseteq \mathbb{R }^{n+1} $ for some $r > 0$,
   and $ \omega \in C^{\bar{n}+1}(\mathcal{B},\mathbb{R})$ satisfy the R\"{u}ssmann's non-degenerate condition with the index $\bar{n}$ and  amount $\beta$. 
   Then for aforementioned  $\omega_v, v=0,1,2,3 \dots$ , 
   let 
   \begin{equation} 
    \bar{V}_v= \big \{ \xi \in V : | \langle k,t\omega_v \rangle -2\pi l| \geq \frac{t\gamma_v}{|k|^{\tau}}   , \forall \  k\in \mathbb{Z}^n \setminus \{0\}, \forall \ l\in \mathbb{Z} \big \}
   \end{equation}
   \begin{equation} \label{201}
       \bar V = \bigcap \limits_{v=0} ^{\infty} \bar{V}_v   
   \end{equation}
   then we have the estimate
   \begin{equation*}
     \big |V \setminus \bar V \big | \leq c_3   \gamma  d^{n}((n+1)^{-\frac{1}{2}}+2d+ \theta^{-1}d)
   \end{equation*}
where $c_3 = 12 \bar{c} (2\pi e)^\frac{n+1}{2} (\bar n +2)^{\bar n +3} [(\bar n  +2)!]^{-1} \big(\lvert \omega \rvert_V + c  \epsilon +1 \big) \beta^{-\frac{\bar{n}+2}{\bar{n } +1}} \lvert  \omega \rvert^{\bar{ n}  +2}_{\mathcal{B}} $, $\bar{c}$ depends only on $n$, $\bar{n}$ and $\tau$.
 And for $\gamma$ small enough ,  $\bar{V}$ can be a set of positive measures.
\end{theorem}

Before proving the theorem, we first refer to the following lemma  in \cite{r3} ).

\begin{lemma} \label{l1}
   (Lemma 4.9 in \cite{r3})
   Let $ K \subseteq \mathbb{R}^n $ be compact with positive diameter $d := \sup_{x,y \in K } |x-y|_2 >0$, 
   define $\mathcal{B}:= {K +\theta }\subseteq \mathbb{R }^ n $ for some $\theta > 0$,
   and $ g \in C^{u_0+1}(\mathcal{B},\mathbb{R})$ be a function with 
\begin{equation} \label{e7}
    \min_{y\in K} \max_{0\leq v\leq {u_0}} |D^v g(y)|\geq \beta
\end{equation}
   for some ${u_0} \in \mathbb{N} $ and $\beta>0$. Then for any $\tilde g\in C^{u_0}(\mathcal{B},\mathbb{R})$ satisfying $|\tilde g - g |^{u_0}_\mathcal{B} \leq \frac{1}{2} \beta$ we have the estimate
\begin{equation}
    \big| \big \{ y\in K : |\tilde g(y)|\leq \epsilon \big \} \big | \leq B d^{n-1}(n^{-\frac{1}{2}}+2d+\theta^{-1}d)(\frac{\epsilon}{\beta})^\frac{1}{u_0}\frac{1}{\beta} \max_{0 < v\leq {u_0}} |D^v g|_\mathcal{B}
\end{equation}
    whenever $0<\epsilon<\frac{\beta}{2u_0 +2}$. Here, $B:=3(2\pi e)^\frac{n}{2} (u_0+1)^{u_0+2} [(u_0 +1)!]^{-1}$.
\end{lemma}
\begin{remark}
   This Lemma is Lemma 4.9 in \cite{r3}, which comes from Theorem 17.1 in \cite{r2}, and the proof is similar to the proof of Theorem 17.1 in \cite{r2}.
\end{remark}

Now, let us focus on the proof of Theorem \ref{t1}.
\paragraph{Proof of Theorem \ref{t1}}
Firstly, for $t \omega$ satisfies the R\"{u}ssmann's non-degenerate condition, 
by Remark \ref{re03}, $ \exists \ \bar{n} = \bar{n}(\omega, V) \in \mathbb{N}$ and $ \beta = \beta (\omega, V) >0 $ such that $\min \limits_{\xi \in V} \max \limits_{0\leq v \leq \bar{n}} \lvert D^v \langle k, t \omega(\xi)\rangle \rvert \geq t \beta \lvert k \rvert$ for $\forall \ k\in \mathbb{Z}^n \setminus \{0\} $.

Let $\widetilde{V} = V \times (0,1)$, $\widetilde{\xi} = (\xi, \xi^{'})$, $\xi^{'} \in (0,1)$, $t \widetilde{\omega} = (t \omega, -2 \pi)$, $\widetilde{k} = (k, l)$, then we have $ rank \big \{\partial^i_{\xi} \widetilde{\omega}(\widetilde{\xi}) :  |i| \leq \bar{n}+1 \big \} = n+1 ,  \   \widetilde{\xi} \in \widetilde{V}$, and we can easily get that $\min \limits_{\widetilde{\xi} \in \widetilde{V}} \max \limits_{0\leq v\leq \bar{n}+1} |D^v \langle \widetilde{k}, t \widetilde{\omega}(\widetilde{\xi}) \rangle |\geq t \beta \lvert \widetilde{k} \rvert$.

Let $\mathcal{B} = (\widetilde{V}+ r) \cap \mathbb{R}^{n+1}$, by  $\lVert  \omega_v -  \omega \rVert ^* \leq c   \epsilon $, then for $\epsilon$  small enough, we have
$ \lvert \langle \widetilde{k}, t \widetilde{\omega}_v \rangle  - \langle \widetilde{k}, t \widetilde{\omega}  \rangle \rvert^{\bar{n}+1}_{\mathcal{B}}  \leq \frac{1}{2} t \beta \lvert \widetilde{k} \rvert
$, and for $\gamma_0$  small enough, we have $t\gamma_v \lvert k\rvert^{-\tau} \leq \frac{t \beta \lvert \widetilde{k} \rvert}{2 (\bar{n}+1) +2}$, so we can use Lemma \ref{l1}, and we have
\begin{align*}  
    & \big \lvert  \{ \xi \in V : \lvert \langle k, t\omega_v \rangle - 2 \pi l   \rvert \leq t\gamma_v \lvert k\rvert^{-\tau}  \} \big \rvert =  \big \lvert  \{ \widetilde{\xi} \in \widetilde{V} : \lvert \langle \widetilde{k}, t \widetilde{\omega}_v \rangle \rvert \leq t\gamma_v \lvert k\rvert^{-\tau}  \} \big \rvert
    \\
    & \ \leq \ c_1 d^{n}((n+1)^{-\frac{1}{2}}+2d+ r^{-1}d) (\frac{t \gamma_v}{t \beta})^\frac{1}{\bar n+1} \lvert k \rvert^{\frac{ -\tau } {\bar{n} +1} } \lvert \widetilde{k} \rvert^{\frac{ -1 } {\bar{n} +1} }  \frac{1}{t \beta} \max_{0 < v\leq \bar{n}+2} |D^v t \widetilde{\omega}|_\mathcal{B}
\end{align*}
for $  k\in \mathbb{Z}^n \setminus \{0\} ,  l \in \mathbb{Z}$ , where $ c_{1} = 3(2\pi e)^\frac{n+1}{2} (\bar n +2)^{\bar n +3} [(\bar n  +2)!]^{-1} $, $d$ is the diameter of $V \times (0,1)$.
For such $l$ right here, we have $\lvert l \rvert \leq \lvert k \rvert  \hat{K} $, where $\hat{K} = \lvert \omega \rvert_{V} + c \epsilon +1$.  
Now we define  
 \begin{align*}
    R^{k,l}_v &= \big \{ \xi \in V :  \lvert \langle k, t\omega_v \rangle -2 \pi l  \rvert  < t\gamma_v \lvert k \rvert^{-\tau}  \big \}, \ \lvert l \rvert \leq \lvert k \rvert  \hat{K}  , \  k\in \mathbb{Z}^n \setminus \{0\}  \\
    R^k_v &= \big \{ \xi \in V :  \lvert \langle k, t\omega_v \rangle -2 \pi l  \rvert  < t\gamma_v \lvert k \rvert^{-\tau} , \exists \ l \in \mathbb{Z}   \big \}, \  k\in \mathbb{Z}^n \setminus \{0\} \\
    R_v &= \big \{ \xi \in V :  \lvert \langle k, t\omega_v \rangle  -2 \pi l  \rvert  < t\gamma_v \lvert k \rvert^{-\tau} ,\exists \ l \in \mathbb{Z}  ,  \exists  \ k\in \mathbb{Z}^n \setminus \{0\} \big \}
\end{align*}
i.e., $ R_v= \bigcup \limits_{k\in \mathbb{Z}^n \setminus \{0\}} R^k_v = \bigcup \limits_{k\in \mathbb{Z}^n \setminus \{0\} } \bigcup \limits_{\lvert l \rvert \leq \lvert k \rvert  \hat{K} }  R^{k,l}_v$.
Besides, for $ k\in \mathbb{Z}^n \setminus \{0\}$ we have $\lvert \widetilde{k} \rvert^{\frac{ -1 } {\bar{n} +1} } \leq 1$ , then,
\begin{align*}
    \lvert R^{k,l}_v \rvert & \leq c_1 d^{n}((n+1)^{-\frac{1}{2}}+2d+ r^{-1}d) (\frac{t \gamma_v}{t \beta})^\frac{1}{\bar n+1}  \frac{1}{t \beta} \max_{0 < v\leq \bar{n}+2} |D^v t \widetilde{\omega}|_\mathcal{B} \lvert k \rvert^{\frac{ -\tau} {\bar n +1 }} \\
    & \leq c_1 d^{n}((n+1)^{-\frac{1}{2}}+2d+ r^{-1}d) (\frac{ \gamma_v}{ \beta})^\frac{1}{\bar n+1}  \frac{1}{ \beta} \max_{0 < v\leq \bar{n}+2} |D^v  \omega|_\mathcal{B} \lvert k \rvert^{\frac{ -\tau} {\bar n +1 }} \\
    \lvert R^k_v \rvert & = \lvert \bigcup \limits_{\lvert l \rvert \leq \lvert k \rvert  \hat{K} }  R^{k,l}_v \rvert  \leq 2 \lvert k \rvert \hat{K} \lvert R^{k,l}_v \rvert \\
    & \leq 2 \hat{K} c_1   d^{n}((n+1)^{-\frac{1}{2}}+2d+ r^{-1}d) (\frac{ \gamma_v}{ \beta})^\frac{1}{\bar {n}+1}  \frac{1}{ \beta} \lvert  \omega \rvert^{\bar{ n}  +2}_{\mathcal{B}} \lvert k \rvert^{\frac{\bar{n} +1 -\tau} {\bar {n} +1 }} \\
    & \leq c_2   d^{n}((n+1)^{-\frac{1}{2}}+2d+ r^{-1}d) (\frac{ \gamma_v}{ \beta})^\frac{1}{\bar {n}+1}  \frac{1}{ \beta} \lvert  \omega \rvert^{\bar{ n}  +2}_{\mathcal{B}}  \lvert k \rvert^{\frac{\bar{n} +1 -\tau} {\bar {n} +1 }} 
\end{align*}
where $c_2 = 6 (2\pi e)^\frac{n+1}{2} (\bar n +2)^{\bar n +3} [(\bar n  +2)!]^{-1} \hat{K}$.

By (\ref{201}) $ \bar V = \bigcap \limits_{v=0} ^{\infty} \bar{V}_v $,
we have the following inequality.
\begin{align*}
    & \lvert V \setminus \bar{V} \rvert  \leq \sum _{v=0}^\infty \lvert R_v \rvert  \ \leq \sum _{v=0}^\infty  \sum _{k\in \mathbb{Z}^n \setminus \{0\} } \lvert R_v^k \rvert 
    \leq \sum _{v=0}^\infty \sum^\infty _{r=1} 2n(2r+1)^{n-1}\lvert R_v^r \rvert \\
   & \leq \sum _{v=0}^\infty \sum^\infty _{r=1} 2n(2r+1)^{n-1} c_2   d^{n}((n+1)^{-\frac{1}{2}}+2d+ \theta^{-1}d) (\frac{ \gamma_v}{ \beta})^\frac{1}{\bar {n}+1}  \frac{1}{ \beta} \lvert  \omega \rvert^{\bar{ n}  +2}_{\mathcal{B}}  r^{\frac{\bar{n} +1 -\tau} {\bar {n} +1 }}   \\
   &\leq \sum ^\infty _{r=1} 2n(2r+1)^{n-1} r^{\frac{\bar{n} +1 -\tau} {\bar {n} +1 }}  c_2   d^{n}((n+1)^{-\frac{1}{2}}+2d+ \theta^{-1}d) \\
   & \ \cdot \sum _{v=0}^\infty \gamma_v^{\frac{1}{\bar {n}+1}}  \beta^{-\frac{\bar{n}+2}{\bar{n } +1}} \lvert  \omega \rvert^{\bar{ n}  +2}_{\mathcal{B}} 
\end{align*}
As for $\tau \geq (n+2)(\bar{n}+1)$, $\bar{c} :=  \sum \limits_{r=1}^\infty  2n(2r+1)^{n-1} r^{\frac{\bar{n} +1 -\tau} {\bar {n} +1 }} $ is convergent. 
Besides, $\gamma_v \leq  \frac{\gamma ^{\bar{n}+1}}{2^{(\bar n +1)v}} $, 
we have $ \sum \limits_{v=0}^\infty \gamma_v^{\frac{1}{\bar n+1}} \leq \sum \limits_{v=0}^\infty \frac{\gamma}{2^v} = 2 \gamma $, and $\lvert  \omega \rvert^{\bar{ n}  +2}_{\mathcal{B}}  < \infty$ by (\ref{001}). Therefore, we have
\begin{align*}
    \lvert V \setminus \bar{V} \rvert & \leq \bar{c} \cdot 2\gamma \cdot   c_2  d^{n}((n+1)^{-\frac{1}{2}}+2d+ \theta^{-1}d)  \beta^{-\frac{\bar{n}+2}{\bar{n } +1}} \lvert  \omega \rvert^{\bar{ n}  +2}_{\mathcal{B}}  \\
    & \leq c_3   \gamma  d^{n}((n+1)^{-\frac{1}{2}}+2d+ \theta^{-1}d)
\end{align*}
where $c_3 = 12 \bar{c} (2\pi e)^\frac{n+1}{2} (\bar n +2)^{\bar n +3} [(\bar n  +2)!]^{-1} \big(\lvert \omega \rvert_V + c  \epsilon +1 \big) \beta^{-\frac{\bar{n}+2}{\bar{n } +1}} \lvert  \omega \rvert^{\bar{ n}  +2}_{\mathcal{B}} $, $\bar{c}$ depends only on $n$, $\bar{n}$ and $\tau$. 
So, for $\gamma$  small enough ,  $\bar{V}$ can be a set of positive measures. 

The proof of Theorem \ref{t1} is complete. $\square$

So we get the measure estimate that  $t \omega_v, v=0,1,2,3 \dots$  satisfy (\ref{e1}). 

Now, let us focus on (\ref{e2}). By the way, (\ref{e3}) and (\ref{e4}) are similar to that. 
Firstly, similar to (\ref{e1}), we focus on
\begin{equation*}
     \big \{\xi \in V : \lvert  \langle k, t \omega_v(\xi) \rangle - \theta_j^{t,v} -2 \pi l   \rvert \geq \frac{t \gamma_v}{{\lvert k \rvert}^\tau} , \forall \  k \in \mathbb{Z}^n \setminus \{0 \}, \forall  \ l \in \mathbb{Z} \big \}
\end{equation*}
We estimate its measure, and we give the following theorem.

\begin{theorem} \label{t2}
       For the same notations in Theorem \ref{t1}, and $\theta_j^{t,v}(\xi), j=1,2, \dots, m$ are the ones above, 
   let 
   \begin{equation*}
        \bar{\bar{V}}_{v} = \big \{ \xi \in V : | \langle k,t\omega_v \rangle - \theta_j^{t,v} -2\pi l| \geq \frac{t\gamma_v}{|k|^{\tau}}  , \forall  k \in \mathbb{Z}^n \setminus \{0 \}, \forall l\in \mathbb{Z}, \forall 1 \leq j \leq m \big \} 
   \end{equation*}
   \begin{equation} \label{202}
    \bar{\bar{V}} = \bigcap \limits_{v=1} ^{\infty} \bar{\bar{V}}_{v} 
    \end{equation}
   then we have the estimate
\begin{equation}
      \big |V \setminus \bar{\bar{V}} \big | 
      \leq  c_5   \gamma  d^{n} 
\end{equation}
here $c_5 = 4  \bar{\bar{c}}  c_4 m (\lvert \omega \rvert_{V} + c  \epsilon + \max \limits_{1 \leq j \leq m}3|B^0_j |  +1) $, $\bar{\bar{c}}$ depends only on $n$, $\bar{n}$ and $\tau$ , $c_4$ and  $c$  are independent of $t$, $\epsilon$ and $\gamma$. And for $\gamma$  small enough ,  $\bar{\bar{V}}$ can be a set of positive measures.
\end{theorem}

\paragraph{Proof of Theorem \ref{t2}} \label{p2}

Review the previous, $ \underline{A_t} = I +  t(A -I), B_t =t B $. We need adjust them to satisfy ${\underline{A_t}}^2- {B_t}^2=I$.
Here $\underline{A_{t j}} = \sec \theta^t _j $, $B_{t j} =\tan \theta^t _j$ for $j = 1,2, \dots , m$. When $t \to 0$, we have $\theta^t_j (\xi) \to 0$, and it can be controlled by $3 |B^0_j t |$. 
This is due to
\begin{equation}
    \lim_{t \to 0} \frac{| \theta_j ^t|}{t}=\lim_ {t\to0} \frac {\arctan{(| t B_j |)}}{t}=\lim_{t \to 0}\frac{|B_j |}{t^2 {B_j}^2 +1}= |B_j |
\end{equation}
So for $t$ small enough, we have $| \theta_j ^t | \leq \frac{3}{2} |B_j t |$, namely, $| \theta_j ^{t,v} | \leq \frac{3}{2} |B_j^v t |$.  
And from the Remark \ref{re10} ,  $ | B^v t | \leq 2 | B^0 t|$, namely, $ | B^v_j t| \leq 2 | B^0_j t|$ .
Thus $| \theta^{t, v} _j| \leq 3|B^0_j t|$, namely,
$|\theta ^{t, v} | \leq 3 |B^0 t|$. 
Then, to complete the proof, we give some lemmas.

\begin{lemma} \label{le001}
    Suppose function $g(x)$ is  real analytic on the closure $\bar{I}$ of $I$, where $I \subset \mathbb{R}$ is an internal. $\lvert g^{(m)}(x) \rvert \geq d$ on $I$, where $d>0$ is a constant. Let $I_h = \{ x \in I : \lvert g(x) \rvert < h$, then we have $m e a s (I_h) \leq 2 (\frac{m! h}{d})^{\frac{1}{m}}$.
\end{lemma}
\paragraph{Proof} 
Note that $I_h$ has at most one connected component, we can assume that $\bar{I_h} = [x_0 -  \sigma_0, x_0 +  \sigma_0]$ with $\sigma_0 \geq 0$. Then, if $\sigma_0 >0$, we have $d(x_0, \partial I_h) > \sigma$ for $\forall \ 0< \sigma < \sigma_0$, then on $I_h$, by Cauchy inequality, we have $d \leq \lvert g^{m}(x) \rvert \leq \frac{m! h}{\sigma^m}$. By it holds for $\forall \ 0< \sigma < \sigma_0$, so $\lvert I_h \rvert \leq 2 \sigma_0 \leq 2 (\frac{m! h}{d})^{\frac{1}{m}}$. If $\sigma_0 = 0$, the result still holds. $\square$

Next, refer to Theorem 1.2 in \cite{r4}  given by Zhang D, Wang S and Xu J, we get the following lemma.

\begin{lemma} \label{l2}
    Suppose $\omega(\xi)$ is real analytic on $\Sigma$ and satisfy the R\"{u}ssmann's non-degenerate condition: $\langle \omega , y \rangle \not\equiv 0 $ for all $y \in \mathbb{R}^n \setminus \{0\}$ .  Let $P(\xi)=(P_{i j}(\xi))_{N\times N}$ be an $N\times N$ matrix, and let $M$ be a matrix with
   \begin{equation}
       M(\xi)= \ \langle k, t \omega(\xi) \rangle I_N + d i a g (\theta^t_1(\xi),\theta^t_2(\xi),\dots,\theta^t_N(\xi))+P_t(\xi)
   \end{equation}
   where $I_N$ denotes unit matrix, $0< t < 1$. Let
   \begin{equation}
    R_k^t(\alpha)= \big \{ \xi \in \Sigma : \lVert M(\xi)^{-1} \rVert > \frac{|k|^\tau}{t \alpha} \big \}
   \end{equation}
   where $\lVert \cdot \rVert$ denotes the sup-norm of matrix.
   Then, there exist sufficiently large integers $L>0$ and $K_0>0$ such that, if $\theta^t_j(\xi)$ is $L$-th  differentiable with $\lVert \theta^t_j(\xi) \rVert ^L \leq  \sigma (j=1,2,\dots ,N) $ and 
   \begin{equation}
     {\lVert P_t(\xi) \rVert }^L _{\Sigma} =\sup _{\xi \in \Sigma} \max_{1\leq i,j \leq N} {\lVert P^t_{i,j}(\xi) \rVert }^L \leq \epsilon
   \end{equation}
    where $\sigma$ is a positive  constant, then, for sufficiently small $\epsilon >0$ and $\alpha >0$, we have 
    \begin{equation}
        meas(R_k^t(\alpha)) \leq  c  t^{-\frac{N-1}{L}} [diam ( \Sigma)]^{n-1} (\frac{ \alpha}{|k|^{\tau -N}})^{\frac{1}{L}}, \forall \ |k|> K_0
    \end{equation}
    where $\tau > n L+N $ is a fixed constant and $c$ is independent of  $t$, $\epsilon$ and $\alpha$.
\end{lemma}

\begin{remark} \label{re4}
    The proof of this lemma is exactly like the proof of Theorem 1.2 in \cite{r4}, except that Lemma \ref{le001} in this paper is used instead of Lemma 2.2 in \cite{r4}. Besides,  as $c$ is independent of $\epsilon$, we can replace condition ${\lVert P_t(\xi) \rVert }^L _{\Sigma} \leq \epsilon$ with ${\lVert P_t(\xi) \rVert }^L _{\Sigma} \leq \lvert k \rvert \epsilon$, the conclusion still holds for fixed $k$. See \cite{r4} for specific proof. 
\end{remark}

Extend the notation in the proof of Theorem \ref{t1}, e.g., $d$, $\bar{n}$, $\beta$, $\widetilde{V} = V \times (0,1)$, $\widetilde{\xi} = (\xi, \xi^{'})$,  $t \widetilde{\omega} = (t \omega, -2 \pi)$, $\widetilde{k} = (k, l)$. 
Now we let $P_{v,l}^t(\widetilde{\xi}) = \langle \widetilde{k}, t \widetilde{ \omega}_v(\widetilde{\xi}) - t \widetilde{ \omega}(\widetilde{\xi}) \rangle $be an $1 \times 1$ matrix.  
 For $ \forall \ j \in \{1, 2, \dots, m \}$, let
\begin{equation*}
    M_{v,l}^t(\widetilde{\xi}) = \langle \widetilde{k}, t \widetilde{ \omega}(\widetilde{\xi}) \rangle  - \theta_j^{t,v}(\widetilde{\xi}) + P_{v,l}^t(\widetilde{\xi})
\end{equation*}
Then, by Lemma \ref{l2} and Remark \ref{re4}, we know that there are such  $L>0$ and $K_0>0$.  
From the previous analysis, $| \theta^{t, v} _j| \leq 3|B^0_j |t$, it is easy to get $\lVert \theta_j^{t,v} \rVert ^L \leq  \sigma (j=1,2,\dots ,m) $ with $\sigma$ is a positive constant.
And by $\lVert t \omega_v - t \omega \rVert ^* \leq c t  \epsilon $, 
we have ${\lVert P_{v,l}^t(\widetilde{\xi}) \rVert }^L _{\widetilde{V}} \leq c^{'} t |\widetilde{k}| \epsilon$ with $ \epsilon$  small enough. 
Then for $\gamma_v$ small enough, we have 
\begin{equation}
    \lvert \big \{ \widetilde{\xi} \in \widetilde{V} : \lVert M_v^t(\widetilde{\xi})^{-1} \rVert > \frac{|\widetilde{k}|^\tau}{t \gamma_v} \big \} \rvert \leq c_4 d^n (\frac{ \gamma_v}{|\widetilde{k}|^{\tau -1}})^{\frac{1}{L}}, \forall \  |\widetilde{k}| = \lvert k \rvert + \lvert l \rvert > K_0
\end{equation}
where $\tau > (n+1) L + 1$, $d$ is the diameter of $V \times (0,1)$. 
Similar to the proof of Theorem \ref{t1}, we have $\lvert l \rvert \leq \lvert k \rvert  \hat{K}^{'} $ with $\hat{K}^{'} =  \lvert \omega \rvert_{V} + c \epsilon + \max \limits_{1 \leq j \leq m}3|B^0_j | +1$.

Define  
 \begin{align*}
    R^{'}_{k,v,l} &= \big \{ \xi \in V:  \lvert \langle k, t \omega_v \rangle - \theta_j^{t,v} -2 \pi l   \rvert  < \frac{t\gamma_v}{\lvert k \rvert^{\tau}}   \big \}, \ \lvert l \rvert \leq \lvert k \rvert  \hat{K}^{'}  , \  \lvert k \rvert > K_0 \\ 
    R^{'}_{k,v} &= \big \{ \xi \in V :  \lvert \langle k, t\omega_v \rangle -  \theta_j^{t,v} -2 \pi l  \rvert  < \frac{t\gamma_v}{\lvert k \rvert^{\tau}} , \exists \ l \in \mathbb{Z}   \big \}, \  \lvert k \rvert > K_0 \\
    R^{'}_{v} &= \big \{ \xi \in V :  \lvert \langle k, t\omega_v \rangle - \theta_j^{t,v} -2 \pi l  \rvert  < \frac{t\gamma_v}{\lvert k \rvert^{\tau}} ,\exists \ l \in \mathbb{Z}  ,  \exists  \ \lvert k \rvert > K_0  \big \}
\end{align*}
i.e., $ R^{'}_{v} = \bigcup \limits_{\lvert k \rvert > K_0 } R^{'}_{k,v} = \bigcup \limits_{\lvert k \rvert > K_0  } \bigcup \limits_{\lvert l \rvert \leq \lvert k \rvert  \hat{K}^{'} }  R^{'}_{k,v,l}$.

Then we have
\begin{equation*}
    \lvert R^{'}_{k,v,l} \rvert  \leq c_4 d^{n} (\frac{ \gamma_v}{|\widetilde{k}|^{\tau -1}})^{\frac{1}{L}} \leq c_4 d^{n}  \gamma_v^{\frac{1}{L}} \lvert k \rvert^{- \frac{ \tau -1 } {L} }
\end{equation*}
By  $ R^{'}_{k,v} =  \bigcup \limits_{\lvert l \rvert \leq \lvert k \rvert  \hat{K} }  R^{'}_{k,v,l} $,
 we have
\begin{equation*}
    \lvert R^{'}_{k,v} \rvert \leq 2 \lvert k \rvert \hat{K}^{'} \lvert R^{'}_{k,v,l} \rvert \leq 2 c_4 \hat{K}^{'}  d^{n}   \gamma_v^{\frac{1}{L}} \lvert k \rvert^{- \frac{ \tau -L -1 } {L} } 
\end{equation*}
Let
\begin{align*}
    \bar{\bar{V}}_{j,v}^{'} &= \big \{ \xi \in V : | \langle k,t\omega_v \rangle - \theta_j^{t,v} -2\pi l| \geq \frac{t\gamma_v}{|k|^{\tau}}  , \forall \ |k| > K_0, \forall \ l\in \mathbb{Z} \big \} \\
    \bar{\bar{V}}^{'}_j &= \bigcap \limits_{v=0} ^{\infty} \bar{\bar{V}}_{j,v}^{'} ,  \  \bar{\bar{V}}^{'} = \bigcap \limits_{j=1} ^{m} \bar{\bar{V}}^{'}_j \\
    \bar{\bar{V}}^{''}_{k,v,l,j} &= \big \{ \xi \in V : | \langle k,t\omega_v \rangle - \theta_j^{t,v} -2\pi l| \geq \frac{t\gamma_v}{|k|^{\tau}} \big \} , \   0< |k| \leq K_0, \ l\in \mathbb{Z}  \\
   \bar{\bar{V}}^{''} &= \bigcap \limits_{ k \in \mathbb{Z}^n \setminus \{0\}} \bigcap \limits_{l \in \mathbb{Z}}   \bigcap \limits_{j=1} ^{m} \bigcap \limits_{v=1} ^{\infty} \bar{\bar{V}}^{''}_{k,v,l,j} 
\end{align*}
then we have the following inequality.
\begin{align*}
     \lvert V \setminus \bar{\bar{V}}^{'}_j \rvert & \leq \sum _{v=0}^\infty \lvert R_v^{'} \rvert  \ \leq \sum _{v=0}^\infty  \sum _{\lvert k \rvert > K_0  } \lvert R_{k,v}^{'} \rvert \\
   & \leq \sum _{v=0}^\infty \sum^\infty _{r=K_0 +1} 2n(2r+1)^{n-1}\lvert R_{r,v}^{'} \rvert \\
   & \leq \sum _{v=0}^\infty \sum^\infty _{r=K_0 +1} 2n(2r+1)^{n-1} \cdot 2 c_4 \hat{K}^{'}  d^{n}   \gamma_v^{\frac{1}{L}} r^{- \frac{ \tau -L -1 } {L} }  \\
   &\leq \sum ^\infty _{r=1} 2n(2r+1)^{n-1} r^{- \frac{ \tau -L -1 } {L} }  2 c_4 \hat{K}^{'}  d^{n}    \sum _{v=0}^\infty \gamma_v^{\frac{1}{L}}  
\end{align*}
For $\tau \geq (n+2)L +1$, $\bar{\bar{c}} := \sum ^\infty _{r=1} 2n(2r+1)^{n-1} r^{- \frac{ \tau -L -1 } {L} }$ is convergent. 
And by $\gamma_v \leq  \frac{\gamma ^{L}}{2^{L v}} $, 
we have $ \sum \limits_{v=0}^\infty \gamma_v^{\frac{1}{L}} \leq \sum \limits_{v=0}^\infty \frac{\gamma}{2^v} = 2 \gamma $. Therefore, we have
\begin{equation*}
    \lvert V \setminus \bar{\bar{V}}^{'}_j \rvert  \leq \bar{\bar{c}} \cdot 2 \gamma \cdot  2 c_4 \hat{K}^{'}  d^{n} 
\end{equation*}
By $j$ is arbitrary and $  \bar{\bar{V}}^{'} = \bigcap \limits_{j=1} ^{m} \bar{\bar{V}}^{'}_j $, then we have 
\begin{equation*}
    \lvert V \setminus \bar{\bar{V}}^{'} \rvert \leq m \cdot 4  \bar{\bar{c}} c_4 \hat{K}^{'} \gamma  d^{n}   \leq c_5   \gamma  d^{n} 
\end{equation*}
where $c_5 = 4  \bar{\bar{c}}  c_4 m (\lvert \omega \rvert_{V} + c  \epsilon + \max \limits_{1 \leq j \leq m}3|B^0_j |  +1) $, $c_4$ comes from Lemma \ref{l2} and it is independent of $t$, $\epsilon$ and $\gamma$. 
So, for $\gamma$  small enough ,  $\bar{V}^{'}$ can be a set of positive measures. 

So we get the measure estimate that $\omega_v$ satisfies the condition (\ref{e2}) for $|k|> K_0$. 

Next, we consider $0 < |k| \leq  K_0$. It is mainly finite thought.

Firstly, for  $0 <|k|\leq K_0$, there are only a finite number of $k$ 's like this.  If $\lvert \langle k, t \omega_v(\xi) \rangle  - \theta_j^{t,v}(\xi) -2 \pi l   \rvert  \geq \alpha >0$ on $\bar{\bar{V}}_*$ for  $ \forall \ 0 <|k|\leq K_0$, $ \forall \ \lvert l \rvert \leq K_0 \hat{K}^{'} $ , $ \forall \ 1 \leq j \leq m $ and $ \forall \ v \in \mathbb{N}$, where $\hat{K}^{'} = \lvert \omega \rvert_{V} + c \epsilon + \max \limits_{1 \leq j \leq m}3|B^0_j | +1$,  then by finiteness, we can choose the smallest $\alpha$ like this. As for $\gamma_v$ small enough, we can get $\alpha > t\gamma _v > \frac{t \gamma_v}{|k|^\tau}$. 
Then $\lvert \langle k, t \omega_v(\xi) \rangle   - \theta_j^{t,v}(\xi) -2 \pi l   \rvert \geq \frac{t \gamma_v}{|k|^\tau}$ is satisfied on $\bar{\bar{V}}^{''}$. 
So under the circumstance, we can get that $\bar{\bar{V}}_* \subseteq \bar{\bar{V}}^{''}$.

Secondly, let us deal with the case where this condition is not met. Namely,  $ \exists \ 0 <|k|\leq K_0$, $ \exists \ \lvert l \rvert \leq K_0 \hat{K}^{'} $, $ \exists \ 1 \leq j \leq m $ and $ \exists \ v \in \mathbb{N}$,  there is no such $\alpha$ that  $\lvert \langle k, t \omega_v(\xi) \rangle - \theta_j^{t,v}(\xi) -2 \pi l   \rvert \geq \alpha $ on $V \setminus \bar{\bar{V}}_*$.
Let $ R^{''}_{k,v,l,j} = \big \{ \xi \in V  :  \lvert \langle k, t \omega_v(\xi) \rangle - \theta_j^{t,v}(\xi) -2 \pi l   \rvert  =0  \big \}$,  
by Lemma 18.6 in \cite{r2}, we have  $ \big \lvert R^{''}_{k,v,l,j} \big \rvert = 0$, so $ \lvert V \setminus \bar{\bar{V}}_* \rvert = 0 $, then we have  $ \lvert V \setminus \bar{\bar{V}}^{''} \rvert =0 $.

This is the measure estimate that $\omega_v$ satisfies the condition (\ref{e2}) for $0 < |k| \leq K_0$.

As $\bar{\bar{V}} = \bar{\bar{V}}^{'} \cap \bar{\bar{V}}^{''} $ , we have
\begin{equation}
       \lvert V \setminus \bar{\bar{V}} \rvert \leq   \lvert V \setminus \bar{\bar{V}}^{'} \rvert + \lvert V \setminus \bar{\bar{V}}^{''} \rvert
      \leq  c_5   \gamma  d^{n} 
\end{equation}
$c_5 = 4  \bar{\bar{c}}  c_4 m (\lvert \omega \rvert_{V} + c  \epsilon + \max \limits_{1 \leq j \leq m}3|B^0_j |  +1) $,  $\bar{\bar{c}}$ depends only on $n$, $\bar{n}$ and $\tau$ ,  $c_4$ and $c$ are independent of $t$, $\epsilon$ and $\gamma$. 
So, for $\gamma$  small enough ,  $\bar{\bar{V}}$ can be a set of positive measures. 
The proof of Theorem \ref{t2} is complete. $\square$

So we get the measure estimate that  $t \omega_v, v=0,1,2,3 \dots$  satisfy (\ref{e2}). 

As for (\ref{e3}) and (\ref{e4}), we can do a similar thing.
By the previous analysis in the proof of Theorem \ref{t2},
we know that $| \theta^{t, v} _j| \leq 3|B^0_j t|$, 
 further to have $\lVert \theta^{t, v}_j (\xi) \rVert ^L \leq  \sigma (j=1,2,\dots ,m, v \in \mathbb{N}) $. 
So we can easily get that $\lVert \theta^{t, v}_i(\xi) \pm \theta^{t, v}_j(\xi) \rVert ^L \leq  \sigma (j=1,2,\dots ,m , v \in \mathbb{N}) $, then (\ref{e3}) and  (\ref{e4}) are naturally satisfied. Then we get the following theorem.
\begin{theorem} \label{t5}
       For the same notations in Theorem \ref{t1} and Theorem \ref{t2}, assume that (\ref{e1}), (\ref{e2}),(\ref{e3}) and (\ref{e4}) satisfied on $V_{\gamma}$ for  $ \forall \  k \in \mathbb{Z}^n \setminus \{0 \}$, $ \forall \  l \in \mathbb{Z} $ , $ \forall \ 1 \leq j \leq m $ and $ \forall \ v \in \mathbb{N}$,
          then we have the estimate
   \begin{equation*}
     \lvert V \setminus V_{\gamma} \rvert \leq c_3   \gamma  d^{n}((n+1)^{-\frac{1}{2}}+2d+ \theta^{-1}d) + c_5   \gamma  d^{n} 
   \end{equation*}
where $c_3 = 12 \bar{c} (2\pi e)^\frac{n+1}{2} (\bar n +2)^{\bar n +3} [(\bar n  +2)!]^{-1} \big(\lvert \omega \rvert_{\mathcal{K}} + c  \epsilon +1 \big) \beta^{-\frac{\bar{n}+2}{\bar{n } +1}} \lvert  \omega \rvert^{\bar{ n}  +2}_{\mathcal{B}} $, $c_5 = 4  \bar{\bar{c}}  c_4 m (\lvert \omega \rvert_{V} + c  \epsilon + \max \limits_{1 \leq j \leq m}3|B^0_j |  +1) $, $\bar{c}$ and $\bar{\bar{c}}$ depend only on $n$, $\bar{n}$ and $\tau$,   $c_4$ and  $c$  are independent of $t$, $\epsilon$ and $\gamma$. And for $\gamma$  small enough ,  $V_{\gamma}$ can be a set of positive measures, and $m e a s(V \setminus V_{\gamma}) \to 0$ as $\gamma \to 0$.
\end{theorem}

As for $\ k\in \mathbb{Z}^n \setminus \{0\}$, we have 
$t\gamma_v |k|^{-\tau} \geq \frac{t\gamma_v}{({1+|k|)}^{\tau}}$, so for $\xi \in V \setminus V_*$,  (\ref{e031}) , (\ref{e032}) ,(\ref{e033}) and (\ref{e034}) are naturally satisfied.

\begin{remark} \label{re11}
As of $k=0$, we consider the small divisors which come from the coefficient matrices. By (\ref{e012}), we just have to consider (\ref{e042}) for $k=0$. For fixed $j$, 
by \cite{r6}, we have
\begin{equation*}
     \det (M^v_k) = 2 \sec \theta^{t, v }_j e^{i \langle k, t \omega_v  \rangle } (\cos{\langle{ k, t \omega_v} \rangle} - \cos{\theta^{t, v }_j})
\end{equation*}
 where $M^v_k$ represents the coefficient matrix of (\ref{e042}) in step $v$. When $k=0$, $|\det (M^v_0)| = 2 (|\sec \theta^{t, v }_j| -1)$.
 By the previous definitions, $\underline{A^v_t} = I +  t(A^v - I)$, $\underline{A_{t, j}^v} = \sec \theta^{t,v} _j $ and $A^v_j = \sec \theta^v_j $, then we have $2 (|\sec \theta^{t, v }_j| -1) = 2 t (|A^v_j|-1) = 2 t (|\sec \theta^v_j|-1) $. 
 By Lemma \ref{re10}, we have $ | B^v_j  | \geq \frac{1}{2} | B^0_j  |$ for $v = 1,2,3 \dots$, where $B^v_j = \tan \theta^v_j $, $B^0_j$ is a constant. 
 So for $t$ and $\gamma_0$ small enough, $v = 1,2,3 \dots$, we have 
 \begin{equation*}
     |\det (M^v_0)| = 2 t (|\sec \theta^v_j|-1) \geq 2 t \{ ( (\frac{1}{2} | B^0_j|)^2 +1)^{\frac{1}{2}} -1 \} \geq \frac{c t^2 \gamma_v^2}{(|0| +1)^{2 \tau}}
 \end{equation*}
which is consistent with  $\ k\in \mathbb{Z}^n \setminus \{0\}$.
And for $v=0$, as $\theta_j^0$ is a constant, then $|\det (M^0_0)| = 2 t (|\sec \theta^0_j|-1) \geq \frac{c t^2 \gamma_0^2}{(|0| +1)^{2 \tau}}$ hold naturally for $t$ and $\gamma_0$ small enough.
Since $j$ is arbitrary, then we have (\ref{e031}) , (\ref{e032}) ,(\ref{e033}) and (\ref{e034}) are satisfied for $k=0$, $ \forall \ 1 \leq j \leq m $ and $ \forall \ v \in \mathbb{N}$.
That is to say, for $k=0$, the small divisor conditions hold naturally.
\end{remark}

Thus, measure estimate in Theorem \ref{t5} is corresponding to all small divisor conditions. The proof of Theorm \ref{t02} is complete. $\square$

Similar to \cite{r6}, the proof of Theorem \ref{t01} is based on Theorem \ref{t02}, and the concrete proof refers to the proof of Theorem 1.1 in \cite{r6}.

\section {Applied to the symplectic difference scheme (\ref{004})}
\label{se6}

Let us focus on the symplectic difference scheme (\ref{004}).

\subsection {Proof of Corollary \ref{co1}}

Now, as for the symplectic difference scheme (\ref{004}), we focus on Corollary \ref{co1} for $\epsilon$ and $t$ small enough.
The case is simple, 
we just substitute $\epsilon H_1 + t^s P_1$ in (\ref{004}) as a perturbation in (\ref{e01}), 
namely, let $P^{'} = \epsilon H_1 + t^s P_1$ . 
Assume that $\sup_{D} |P^{'} | \leq M_1 \epsilon  +  2 n M_2 l_* t^s \leq \epsilon^{'} $ ,
then by Theorem \ref{t02},
we get the corresponding conclusions with $\epsilon^{'}$ rather than $\epsilon$, including the existence of invariant tori and measure estimation.  Then we get Corollary \ref{co1}. $\square$

\subsection {Proof of Corollary \ref{co2}}

As for the  Corollary \ref{co2}, its proof is obvious. $\square$



\end{document}